\documentclass[12pt]{amsart}
\usepackage{amsmath,amssymb,epsfig,color}
\newtheorem{theorem}{Theorem}[section]

\newtheorem{corollary}[theorem]{Corollary}
\newtheorem{definition}[theorem]{Definition}
\newtheorem{lemma}[theorem]{Lemma}
\newtheorem{proposition}[theorem]{Proposition}
\theoremstyle{remark}
\newtheorem{remark}[theorem]{Remark}
\newtheorem{example}[theorem]{Example}

\newcommand{\vanish}[1]{}\parskip=12pt

\newcommand{\0}{\widehat{0}}
\newcommand{\1}{\widehat{1}}
\newcommand{\rank}{\rho}

\newcommand{\QQ}{\widetilde{Q}}
\newcommand{\st}{\operatorname{t}}
\newcommand{\CC}{\operatorname{{\mathcal C}}}
\newcommand{\DD}{\operatorname{{\mathcal D}}}
\newcommand{\rev}{\operatorname{rev}}

\begin{document}
\title{The short toric polynomial} 
\author{G\'abor Hetyei}
\address{Department of Mathematics and Statistics, UNC Charlotte, 
	Charlotte, NC 28223}
\email{ghetyei@uncc.edu}
\begin{abstract}
We introduce the short toric polynomial associated to a graded Eulerian poset.
This polynomial contains the same information as the two toric polynomials
introduced by Stanley, but allows different algebraic manipulations. 
The intertwined recurrence defining Stanley's toric polynomials may be
replaced by a single recurrence, in which the degree of the discarded
terms is independent of the rank. A short toric variant
of the formula by Bayer and Ehrenborg, expressing the toric
$h$-vector in terms of the $cd$-index, may be stated in a
rank-independent form, and it may be shown 
using weighted lattice path enumeration and the reflection
principle. We use our techniques to derive a formula expressing the
toric $h$-vector of a dual simplicial Eulerian poset in terms of its
$f$-vector. This formula implies Gessel's formula 
for the toric $h$-vector of a cube, and may be used to prove that the 
nonnegativity of the toric $h$-vector of a simple polytope is a
consequence of the Generalized Lower Bound Theorem holding for
simplicial polytopes. 
\end{abstract}
\maketitle

\section*{Introduction}

As mathematicians, we often look for a ``magic'' simplification that makes
known results easier to state, and helps us find new results which were
cumbersome to even talk about using the old terminology. 
In the study of Eulerian partially ordered sets such a wonderful
simplification was the introduction of the $cd$-index by
Fine (see \cite{Bayer-Klapper}) allowing to restate the already known
Bayer-Billera formulas~\cite{Bayer-Billera} in a simpler form and 
to formulate Stanley's famous nonnegativity conjecture~\cite{Stanley-flag}  
regarding the $cd$-coefficients of Gorenstein$^*$ posets, shown many
years later by Karu~\cite{Karu}.

The present author believes that a similar ``magic'' moment has yet to
arrive in the study of the toric polynomials $f(P,x)$ and $g(P,x)$
associated to an Eulerian poset $\widehat{P}=P\uplus\{\1\}$ by
Stanley~\cite{Stanley-GH}. Without doubt, these invariants are very
important, linked to deep results in algebraic topology, and yielding
highly nontrivial combinatorial interpretations, whenever such
interpretations were found. However, the defining intertwined recurrence
is difficult to use directly, not only because two sequences of
polynomials need to be defined simultaneously, but also because 
the degree of the terms to be discarded in the process changes all the
time as the rank of the intervals considered changes.

The introduction of the new invariant proposed in this paper is
probably not the desired ``magic simplification'' yet, but it represents
a modest improvement in some cases. The idea on which it is based is
very simple and could also be used beyond the confines of our
current area. As explained in Section~\ref{sec:sym}, there is a bijective way
to associate each  {\em multiplicatively symmetric} polynomial $p(x)$
(having a symmetric array of coefficients) 
to an {\em additively symmetric polynomial $q(x)$} (whose multiset of
zeros is symmetric to the origin) of the same degree, having the same
set of coefficients. For example, the additively symmetric variant of
$1-2x+7x^3-2x^5+x^6$ is $x^6-2x^4+7$. There is no change when we want to
extract the coefficients of the individual polynomials only, but when we
consider a sequence $\{p_n(x)\}_{n\geq 0}$ of multiplicatively symmetric
  polynomials, given by some rule, switching to
  the additively symmetric variant $\{q_n(x)\}_{n\geq 0}$ greatly
  changes the appearance of the rules, making them sometimes easier to
  manipulate. Since multiplicatively symmetric polynomials abound in
  combinatorics, the basic idea presented in Section~\ref{sec:sym} is
  worth trying in many situations, unrelated to our current subject.

The short toric polynomial $\st(P,x)$, associated to a graded Eulerian
poset $\widehat{P}$ is defined in Section~\ref{sec:sth} as the
additively symmetric variant of Stanley's toric polynomial
$f(P,x)$. The intertwined recurrence defining $f(P,x)$ and $g(P,x)$ is
equivalent to a single recurrence for $\st(P,x)$. In this recurrence,
multiplication by negative powers of $x$ occurs and we obtain a
polynomial by discarding all terms of negative degree and also certain
constant terms. It is a tempting thought to use this recurrence to
generalize the short toric polynomial to all ranked posets having a
unique minimum element, even if in the cases of lower Eulerian posets,
``severe loss of information'' may occur, compared to Stanley's
generalization of $f(P,x)$ to such posets. We state and outline the
proof of the short toric variant of Fine's formula (see \cite{Bayer} and
\cite[Theorem 7.14]{Bayer-Ehrenborg}) expressing the toric $h$-vector in terms
of the flag $f$-vector. Using this formula, it is easy to observe that
the generalization of $\st(P,x)$ makes most sense for ranked posets 
with unique minimum element $\0$ such that the reduced Euler
characteristic of the order complex of $P\setminus\{\0\}$ is not zero. 

Arguably the nicest result in this paper is
Theorem~\ref{thm:CCDD} in Section~\ref{sec:t-cd}, expressing the short
toric polynomial associated to a graded Eulerian poset by defining two
linear operators on the vector space of polynomials that need to be
substituted into the reverse of the $cd$-index and applied to the
constant polynomial $1$. The fact that the toric $h$-vector may be
computed by replacing the letters $c$ and $d$ in the reverse of the
$cd$-index by some linear operators and applying the resulting 
linear operator to a specific vector is a direct consequence 
of the famous result by Bayer and Ehrenborg~\cite[Theorem
  4.2]{Bayer-Ehrenborg}, expressing the toric $h$-vector in terms of the
$cd$-index. In applications, the use of the Bayer-Ehrenborg
result may be facilitated by finding a linearly equivalent presentation 
that is easier to manipulate. In this sense our Theorem~\ref{thm:CCDD}
is analogous to Lee's result~\cite[Theorem 5]{Lee}, presenting another 
easily memorizable reformulation of~\cite[Theorem
  4.2]{Bayer-Ehrenborg}. Our Theorem~\ref{thm:CCDD} offers the first
rank-independent substitution rule, making it 
more useful in proofs involving induction on rank.  
Theorem~\ref{thm:Finecd-t}, which is the reason
behind Theorem~\ref{thm:CCDD}, also implies the 
short toric variant of the Bayer-Ehrenborg result~\cite[Theorem
  4.2]{Bayer-Ehrenborg}, and has a proof using weighted lattice path
enumeration and the reflection principle. The idea of using a weighted
lattice path model to interpret Fine's formula is already present in 
the work Bayer and Ehrenborg~\cite[Section 7.4]{Bayer-Ehrenborg} 
where it is used to provide an alternative proof of their
formula~\cite[Theorem 3.1]{Bayer-Ehrenborg}  
expressing the toric $h$-vector in terms of the $ab$-index. 
By finding the $cd$-index via calculating the $ce$-index first (instead
of the $ab$-index), and by using the short toric form, the
applicability of the reflection principle becomes apparent.

The short toric variant of the above cited Bayer-Ehrenborg result
highlights the importance of a sequence of 
polynomials $\{\QQ_n(x)\}_{n\geq 0}$, a variant of the polynomials 
$\{Q_n(x)\}_{n\geq 0}$ already used by Bayer and
Ehrenborg~\cite{Bayer-Ehrenborg}. In Section~\ref{sec:bases} we take a
closer look at this basis of the vector space of polynomials, alongside
the basis formed by the short toric polynomials $\{t_n(x)\}_{n\geq 0}$
associated to Boolean algebras. The polynomials $\{\QQ_n(x)\}_{n\geq 0}$
turn out to be the dual basis to the {\em Morgan-Voyce polynomials};
whereas the polynomials $\{t_n(x)\}_{n\geq 0}$ may be used to provide a
simple formula connecting the short toric polynomial to
Stanley's toric polynomial $g(P,x)$. 

The proof is always in the pudding; the usefulness of a technique 
is much more apparent if it is used to answer a question that was open
before. Such an application may be found in Section~\ref{sec:duals}
where we express the toric $h$-vector of an Eulerian dual
simplicial poset in terms of its $f$-vector. This question was raised by Kalai,
see~\cite{Stanley-GH}. Besides using Theorem~\ref{thm:CCDD}, the proof
of the formula depends on a formula conjectured by Stanley~\cite[Conjecture
  3.1]{Stanley-flag} and shown by the present author~\cite[Theorem
  2]{Hetyei-Andre}, expressing the contribution of the $h$-vector
entries of an Eulerian simplicial poset to its $cd$-index as weights of
certain Andr\'e permutations. The result was not found easily: 
it was conjectured after using Maple to compute many examples, and then
shown by induction. Finding a ``more combinatorial'' 
explanation in the future is desirable; the numbers indicate that the
models involving Narayana numbers are good candidates for
generalization. Finally, an equivalent form of our formula shows the
following, perhaps surprising result: the nonnegativity of the toric
$h$-vector of simple polytope is a direct consequence of the Generalized
Lower Bound Theorem (GLBT) holding for simplicial polytopes. Thanks to
Karu~\cite{Karu-nr}, we know that the GLBT holds for all polytopes, which
yields a much stronger statement on the $h$-vector of an arbitrary
simple polytope. However, the validity of the GLBT for simplicial
polytopes was shown much earlier by Stanley~\cite{Stanley-NF}; and
the nonnegativity of the toric $h$-entries of a simple polytope 
is derived from this earlier result using a short and elementary
reasoning. 

The paper re-emphasizes the close relation between the study of the
toric polynomials and ``Catalan combinatorics''. This connection 
is already present in the formulas of Bayer and Ehrenborg~\cite[Theorems
  4.1 and 4.2]{Bayer-Ehrenborg} 
expressing the toric $h$-polynomial of an Eulerian poset in terms of its
flag $h$-vector and $cd$-index; the same holds
for the formula of Billera and Brenti~\cite[Theorem 3.3]{Billera-Brenti}
expressing the Kazhdan-Lusztig polynomial of any Bruhat interval in any Coxeter
group in terms of the complete $cd$-index as well as for the formulas of the
present author~\cite{Hetyei-cubical}, expressing the toric
$h$-contributions of cubical shelling components. The present work adds
explicit relations to weighted lattice path enumeration, Morgan-Voyce
polynomials and the Narayana numbers. Besides exploring these
connections further, the generalization of the definition of the short
toric polynomials to non-Eulerian posets is worth further investigation.
A good starting point could be computing the short toric polynomial
associated to the face lattice of a finite dimensional vector space of a
finite field; if the outcome is unsatisfying, the example is probably a
good source of inspiration to define a $q$-analogue.

\section*{Acknowledgments}

I wish to thank Margaret Bayer, Louis Billera and Richard Stanley for
useful advice and encouragement. Christian Krattenthaler's explanation 
greatly increased my understanding of a combinatorial interpretation of
most numbers appearing in my formulas for the toric polynomials of a
dual simplicial poset. Finally, my heartfelt thanks go to an anonymous
referee for the very careful reading of this manuscript and for many
vital corrections.  

\section{Preliminaries}

A partially ordered set $P$ is {\em graded} if it has a unique minimum
  element $\0$, a unique maximum element $\1$ and a {\em rank function
  $\rank: P\rightarrow {\mathbb N}$} satisfying $\rank(\0)=0$ and
  $\rank(y)=\rank(x)+1$ whenever $y$ covers $x$. The rank of $P$
  is $\rank(\1)$. An often studied
  invariant of a graded poset of rank $n+1$ is its {\em flag $f$-vector
    $(f_S\::\: S\subseteq [1,n])$} where $f_S=f_S(P)$ is the number of maximal
    chains in the {\em $S$ rank-selected subposet $P_S=\{u\in P\::\:
      \rank(u)\in S\}$}. Here and throughout this paper we use the interval
    notation $[i,j]$ to denote the set of integers $\{i, i+1,\ldots,
    j\}$. A graded poset is {\em Eulerian} if every
    interval $[u,v]\subseteq P$ with $u<v$ satisfies $\sum_{z\in [u,v]}
    (-1)^{\rank(z)}=0$. All linear relations satisfied by the flag
    $f$-vector of an Eulerian poset were given by Bayer and
    Billera~\cite{Bayer-Billera}. It was observed by Fine and proven
    by Bayer and Klapper~\cite{Bayer-Klapper} that the Bayer-Billera
    relations may be restated as the existence of the {\em $cd$-index},
    as follows. Introducing the {\em flag $h$-vector $(h_S\::\:
      S\subseteq [1,n])$} of a graded poset of rank $(n+1)$ by setting 
\begin{equation}
\label{eq:h-f}
h_S:=\sum_{T\subseteq S} (-1)^{|S|-|T|} f_T,
\end{equation}
  we define the {\em $ab$-index} as the polynomial 
$$
\Psi_{P}(a,b)=\sum_{S\subseteq [1,n]} h_S u_S
$$
in noncommuting variables $a$ and $b$ where the monomial 
$u_S=u_1\cdots u_n$ is given by
$$
u_i=
\left\{
\begin{array}{ll}
b&\mbox{if $i\in S$,}\\
a&\mbox{if $i\not\in S$.}\\
\end{array}
\right.
$$ 
The $ab$-index of an Eulerian poset is then a polynomial of $c=a+b$ and
$d=ab+ba$.  This polynomial $\Phi_P(c,d)$ is called the {\em $cd$-index}
of $P$.  As it was observed by Stanley~\cite{Stanley-flag}, the
existence of the $cd$ index is equivalent to stating that the {\em
  $ce$-index}, obtained by rewriting the $ab$-index as a
polynomial of $c=a+b$ and $e=a-b$, is a polynomial of $c$ and $e^2$. Let
us denote by $L_S$ the coefficient of the $ce$ word $v_1\cdots v_n$ 
given by
$$
v_i=
\left\{
\begin{array}{ll}
e&\mbox{if $i\in S$,}\\
c&\mbox{if $i\not\in S$}\\
\end{array}
\right.
$$   
in the $ce$-index. It was shown by Bayer and
Hetyei~\cite{Bayer-Hetyei-fEuler} that the resulting {\em flag
  $L$-vector $(L_S\::\: S\subseteq [1,n])$} of a graded poset of rank
  $(n+1)$ is connected to the flag $f$-vector by the formulas 
\begin{equation}
\label{eq:fL}
L_S=(-1)^{n-|S|}\sum_{T\supseteq [1,n]\setminus S}
\left(-\frac{1}{2}\right)^{|T|} f_T\quad\mbox{and}\quad
f_S=2^{|S|}\sum_{T\subseteq [1,n]\setminus S} L_T.
\end{equation}
The Bayer-Billera relations are thus also equivalent to stating that, for
an Eulerian poset, $L_S=0$ unless $S$ is an {\em even set}, i.e., a
disjoint union of intervals of even cardinality. A short direct proof 
of this equivalence may be found in~\cite{Bayer-Hetyei-gMoebius}. 

The {\em toric $h$-vector} associated to a graded Eulerian poset
$\widehat{P}=[\0,\1]$ was 
defined by Stanley~\cite{Stanley-GH} by introducing the polynomials 
$f([\0,\1),x)$ and $g([\0,\1),x)$ via the intertwined recurrences 
\begin{equation}
\label{eq:fgrec}
f([\0,\1),x)=\sum_{p\in [\0,\1)} g([\0,p),x) (x-1)^{\rank(\1)-1-\rank(p)}
\quad\mbox{and}
\end{equation}
\begin{equation}
\label{eq:gfrec}
g([\0,\1),x)=\sum_{i=0}^{\lfloor (\rank(\1)-1)/2\rfloor} ([x^i]
  f([\0,\1),x)- [x^{i-1}] f([\0,\1),x)) x^i
\end{equation}
subject to the initial conditions $f(\emptyset,x)=g(\emptyset,x)=1$. 
Here the intervals $P=[\0,\1)$ and $[\0,p)$ are half open: they contain
    the minimum element $\0$ but they do not contain their maximum
    element. In particular, $f(\emptyset,x)$ and
    $g(\emptyset,x)$ are associated to the only Eulerian poset of rank
    $0$. The toric $h$-vector associated to $[\0,\1)$ is then defined as 
the vector of coefficients of the polynomial $f([\0,\1),x)$:
$$
\sum_i h_i x^i:=x^{\rank(\1)-1}f([\0,\1),1/x). 
$$
By \cite[Theorem 2.4]{Stanley-GH} the polynomial above also equals
$f([\0,\1),x)$, the apparently more complicated definition given by
  Stanley~\cite{Stanley-GH} is made for the sake of generalizations to {\em
    lower Eulerian posets}.   
The first formula expressing the polynomials in terms of the flag
$f$-vector was found by Fine (see \cite{Bayer} and \cite[Theorem
  7.14]{Bayer-Ehrenborg}). Here we state it in an equivalent form that
appears in the paper of Bayer and Ehrenborg~\cite[Section 7]{Bayer-Ehrenborg}:
\begin{equation}
\label{eq:Fine}
f([\0,\1),x)=\sum_{S\subseteq [1,n]} f_S \sum_{\lambda\in
    \{-1,1\}^n: S(\lambda)\supseteq S} (-1)^{|S|+n-i_{\lambda}}
  x^{i_{\lambda}},  
\end{equation}
where $S(\lambda)=\{s\in [1,n]\: :\:
\lambda_1+\cdots+\lambda_s>0\}$ and $i_{\lambda}$ is the
number of $-1$'s in $\lambda$. Bayer and Ehrenborg~\cite[Theorem
  4.2]{Bayer-Ehrenborg} also expressed the toric $h$-vector of an
Eulerian poset in terms of its $cd$-index. 

\section{Additive and multiplicative symmetry of polynomials}
\label{sec:sym}

\begin{definition}
Let $p(x)\in K[x]$ be a polynomial of degree
$n$, with coefficients from a field $K$. We say that $p(x)$ is {\em
  multiplicatively symmetric} if it satisfies $x^np(x^{-1})=p(x)$ and
$p(x)$ is {\em additively symmetric} if it satisfies $p(x)=(-1)^n
p(-x)$. 
\end{definition}

\begin{lemma}
A polynomial $p(x)=a_n x^n+a_{n-1}x^{n-1}+\cdots+a_0$ of degree $n$ is
multiplicatively symmetric if and only if its coefficients satisfy
$a_k=a_{n-k}$ for $0\leq k\leq n$.
\end{lemma}

\begin{lemma}
A polynomial $p(x)=a_n x^n+a_{n-1}x^{n-1}+\cdots+a_0$ of degree $n$ is
additively symmetric if and only if its coefficients satisfy
$a_{n-2k-1}=0$ for $0\leq k\leq \lfloor (n-1)/2\rfloor$.
\end{lemma}

\begin{theorem}
\label{thm:addmultsym}
A polynomial $p(x)\in K[x]$ of degree $n$ is multiplicatively symmetric if and
only if there is an additively symmetric polynomial $q(x)\in K[x]$ of degree $n$
satisfying
\begin{equation}
\label{eq:sym}
p(x)=x^{\frac{n}{2}}
\left(q(\sqrt{x})+q\left(\frac{1}{\sqrt{x}}\right)-q(0)\right). 
\end{equation}
Moreover, $q(x)$ is uniquely determined.
\end{theorem}
\begin{proof}
Assume first $p(x)$ may be written in the form given in \eqref{eq:sym}
using an additively symmetric polynomial 
$q(x)=\sum_{k=0}^{\lfloor n/2\rfloor} a_{n-2k} x^{n-2k}$ of degree
$n$. Then we have  
$$
p(x)=\sum_{k=0}^{\lfloor \frac{n}{2}\rfloor} a_{n-2k}
(x^{n-k}+x^k)-x^{\frac{n}{2}} q(0).
$$
For odd $n$ we have  $q(0)=0$, thus $p(x)$ is a polynomial of degree
$n$. Clearly, it is multiplicatively symmetric.

Assume now that $p(x)\in K[x]$ is multiplicatively symmetric of
degree $n$. If $n$ is even then $p(x)$ is of the form 
$$
p(x)=a_{\frac{n}{2}} x^{\frac{n}{2}}+\sum_{k=1}^{\frac{n}{2}} 
a_{\frac{n}{2}-k}\left(x^{\frac{n}{2}+k}+x^{\frac{n}{2}-k}\right)
= x^{\frac{n}{2}}\left(
a_{\frac{n}{2}} +\sum_{k=1}^{\frac{n}{2}} 
a_{\frac{n}{2}-k}\left((\sqrt{x})^{2k}+(\sqrt{x})^{-2k}\right)\right),
$$
and it satisfies \eqref{eq:sym} with 
$q(x)=a_{n/2} +\sum_{k=1}^{n/2} a_{n/2-k} x^{2k}$.

If $n$ is odd then 
$p(x)$ is of the form 
$$
p(x)=\sum_{k=0}^{\frac{n-1}{2}} 
a_{\frac{n-2k-1}{2}}\left(x^{\frac{n-2k-1}{2}}+x^{\frac{n+2k+1}{2}}\right)
= x^{\frac{n}{2}}\left(
\sum_{k=0}^{\frac{n-1}{2}} 
a_{\frac{n-2k-1}{2}}\left((\sqrt{x})^{2k+2}+(\sqrt{x})^{-2k-2}\right)\right),
$$
and it satisfies \eqref{eq:sym} with 
$q(x)=\sum_{k=0}^{(n-1)/2} a_{(n-2k-1)/2} x^{2k+1}$.

We are left to show that $q(x):=\sum_{k=0}^{\lfloor n/2\rfloor} b_k
x^{n-2k}$ is uniquely determined by \eqref{eq:sym}. 
For $k<n/2$, $b_k$ must equal  $[x^{n-k}]p(x)$. Finally, for even
$n$, $b_{n/2}$ must equal $[x^{n/2}] p(x)$.
\end{proof}

\begin{definition}
\label{def:amsym}
Given a multiplicatively symmetric polynomial $p(x)$ we call the {\em
  additively symmetric variant of $p(x)$} the additively symmetric
polynomial $q(x)$ associated to $p(x)$ via \eqref{eq:sym}. Conversely, given   
an additively symmetric polynomial $q(x)$ we call the {\em
  multiplicatively symmetric variant of $q(x)$} the polynomial $p(x)$
defined by \eqref{eq:sym}. 
\end{definition}

It should be noted that the additively symmetric variant of a
multiplicatively symmetric polynomial has the ``same coefficients'',
without the redundant repetitions. For example, the additively symmetric
variant of $p(x)=1-2x^2-2x^5+x^7$ is $q(x)=x^7-2x^3$, and the 
additively symmetric variant of $p(x)=1+2x+x^2$ is $q(x)=x^2+2$.

To express the additively symmetric variant of a multiplicatively
symmetric polynomial, we may use the following {\em truncation operators}.
\begin{definition}
Let $K$ be a fixed field. 
For any $z\in {\mathbb Z}$, the {\em truncation operator}
$U_{\geq z}: K[x,x^{-1}]\rightarrow
K[x,x^{-1}]$ is the linear operator defined on the vector space
$K[x,x^{-1}]$ of Laurent polynomials by discarding all terms of
degree less than $z$. Similarly $U_{\leq z}: K[x,x^{-1}]\rightarrow
K[x,x^{-1}]$ is the linear operator defined by discarding all terms of
degree more than $z$.
\end{definition}
The notation $U_{\geq z}$ and $U_{\leq z}$ is consistent with the
notation used by Bayer and Ehrenborg~\cite{Bayer-Ehrenborg}, who rewrite
\eqref{eq:gfrec} as  
\begin{equation}
\label{eq:g-f}
g([\0,\1),x)=U_{\leq \lfloor n/2 \rfloor} ((1-x)f([\0,\1),x)).
\end{equation} 
\begin{lemma}
\label{l:addmult}
Let $p(x)$ be a multiplicatively symmetric polynomial of degree
$n$. Then the additively symmetric variant $q(x)$ of $p(x)$ satisfies
$$
q(x)=U_{\geq 0}(x^{-n} p(x^2))=U_{\geq 0}(x^n p(x^{-2})).
$$
\end{lemma}
\begin{proof}
Substituting $x^2$ into \eqref{eq:sym} and dividing both sides by $x^n $ yields
$$
x^{-n} p(x^2)=
q(x)+q\left(\frac{1}{x}\right)-q(0).
$$
Applying $U_{\geq 0}$ to both sides results in $q(x)$ on the right hand
side. Finally, $x^{-n} p(x^2)=x^n p(x^{-2})$ is an immediate
consequence of the multiplicative symmetry of $p(x)$.
\end{proof}

\section{The short toric polynomial of an arbitrary graded poset}
\label{sec:sth}

Stanley's generalization~\cite[Theorem 2.4]{Stanley-GH} of the
Dehn-Sommerville equations to the toric polynomial
$f([\0,\1),x)$ associated to an Eulerian poset $[\0,\1]$ states the
  following. 
\begin{theorem}[Stanley]
\label{thm:Stanley-multsym}
For an Eulerian poset $[\0,\1]$ of rank $n+1$, the polynomial $f([\0,\1),x)$ is
multiplicatively symmetric of degree $n$. 
\end{theorem}
\begin{definition}
\label{def:st}
The short toric polynomial $\st([\0,\1),x)$ associated to an Eulerian
  poset $[\0,\1]$ is the additively symmetric variant of the toric
  polynomial $f([\0,\1),x)$.  
\end{definition}
As in the case of $f([\0,\1),x)$, we consider the interval
  $[\0,\1)$ half open. In particular, $f(\emptyset,x)=1$ is restated as 
$\st(\emptyset, x)=1$. 
As an immediate consequence of Definitions~\ref{def:amsym} and \ref{def:st}
we obtain that any Eulerian poset $[\0,\1]$ of rank $n+1$ satisfies 
\begin{equation}
\label{eq:f-t}
f([\0,\1),x)=x^{\frac{n}{2}}
\left(\st([\0,\1),\sqrt{x})+\st\left([\0,\1),\frac{1}{\sqrt{x}}\right)
-\st([\0,\1),0)\right).
\end{equation}
Furthermore, by Lemma~\ref{l:addmult}, we have
\begin{equation}
\label{eq:t-f}
\st([\0,\1),x)=U_{\geq 0}(x^{-n} f([\0,\1),x^2))=U_{\geq 0}(x^n
    f([\0,\1),x^{-2})). 
\end{equation}
\begin{lemma}
\label{l:g-t}
If $[\0,\1]$ is an Eulerian poset of rank $n+1$ then we have 
$$
U_{\geq 1}\left(\st ([\0,\1),x)\cdot\left(x-\frac{1}{x}\right) \right)
=x^{n+1} g\left([\0,\1), x^{-2}\right).
$$
\end{lemma}
\begin{proof}
After substituting $x^{-2}$ into \eqref{eq:g-f}, multiplying both sides by 
$x^{n+1}$, and using the obvious identity 
$x^{n+1} U_{\geq -n}(p(x))=U_{\geq 1}(x^{n+1}p(x))$ we obtain 
$$
x^{n+1}g([\0,\1),x^{-2})=U_{\geq 1} ((x-x^{-1}) x^n f([\0,\1),x^{-2})).
$$
Only terms of degree at least $0$ in $x^n f([\0,\1),x^{-2})$ yield a term
  of degree at least $1$ in $(x-x^{-1}) x^n f([\0,\1),x^{-2})$. Thus we
    may also write 
$$
x^{n+1}g([\0,\1),x^{-2})
=U_{\geq 1} ((x-x^{-1}) U_{\geq 0}(x^n f([\0,\1),x^{-2}))),
$$
and the statement follows from \eqref{eq:t-f}.
\end{proof}
Using \eqref{eq:t-f} and Lemma~\ref{l:g-t} we may show the following
{\em fundamental recurrence} for the short toric polynomial. 
\begin{theorem}
\label{thm:fundrec}
The short toric polynomial satisfies the recurrence
$$
\st([\0,\1),x)=U_{\geq 0}\left((x^{-1}-x)^{\rank(\1)-1}
\!\!+\!\!\!\!\sum_{\0<p<\1} U_{\geq 1}\left(\st ([\0,p),x)
  (x-x^{-1})\right)(x^{-1}-x)^{\rank(\1)-\rank(p)-1}\right).
$$
\end{theorem}
\begin{proof}
Let us set $n:=\rank(\1)-1$. Substituting $x^{-2}$ into \eqref{eq:fgrec}
and multiplying both sides by $x^{n}$ yields  
\begin{equation}
\label{eq:ftrec}
x^{n}f([\0,\1),x^{-2})=\sum_{p\in [\0,\1)} x^{\rank(p)}
    g([\0,p),x^{-2}) 
      (x^{-1}-x)^{\rank(\1)-1-\rank(p)}.
\end{equation}
By \eqref{eq:t-f}, applying $U_{\geq 0}$ to the left hand side yields 
$\st([\0,\1),x)$. We only need to show that the right hand side of
  \eqref{eq:ftrec} is the argument of the operator $U_{\geq 0}$ on the
  right hand side of the statement. Since $g(\emptyset,x)=1$,
 the term associated to $p=\0$ is $(x^{-1}-x)^n$ on the right hand side of
  \eqref{eq:ftrec}. For $p\in (\0,\1)$, the term $x^{\rank(p)}
    g([0,p),x^{-2})$ equals $U_{\geq 1}(\st
      ([\0,p),x)\cdot(x-x^{-1}))$ by Lemma~\ref{l:g-t}.
\end{proof}

\begin{definition}
\label{def:stgen}
We extend the definition of $\st([\0,\1),x)$ to all finite posets $P$ 
that have a unique minimum element $\0$ and a
rank function $\rank: P\rightarrow {\mathbb N}$, satisfying
$\rank(\0)=0$, as follows. 
We set $\st(\emptyset,x):=1$ and 
$$
\st(P,x)=U_{\geq 0}\left((x^{-1}-x)^{n}
\!\!+\!\!\sum_{p\in P\setminus \{\0\}} U_{\geq 1}\left(\st ([\0,p),x)
  (x-x^{-1})\right)(x^{-1}-x)^{n-\rank(p)}\right).
$$
Here $n=\max\{\rank (p)\::\: p\in P\}$.
\end{definition}
Stanley~\cite[\S 4]{Stanley-GH} calls a finite poset {\em lower
  Eulerian} if it has a unique minimum element $\0$ and, for each $p\in
P$, the interval $[\0,p]$ is an Eulerian poset. He extends the
definition of the toric polynomial $f(P,x)$ to lower Eulerian posets by
setting 
$$
f(P,x)=\sum_{p\in P} g([\0,p),x) (x-1)^{n-\rank(p)}.
$$
Here $n$ is the length of the longest chain in $P$. 
\begin{proposition}
\label{prop:lowerE}
Let $P$ be a lower Eulerian poset and let $n$ be the length of the
longest chain in $P$. Then we have
$$
\st(P,x)=U_{\geq 0}(x^n f(P,x^{-2})).
$$  
\end{proposition}
Indeed, when we use the recurrence given in Definition~\ref{def:stgen}
to compute $\st(P,x)$, all intervals $[\0,p)$ on the right hand are
  Eulerian posets with their maximum element removed. Thus we may repeat
  the part of the proof of Theorem~\ref{thm:fundrec} showing that the
  right hand side is $U_{\geq 0}(x^n f(P,x^{-2}))$. 
\begin{remark}
\label{rem:lowerE}
Unlike in the Eulerian case, 
the short toric polynomial $\st(P,x)$ of a lower Eulerian poset $P$ 
may not contain sufficient information to recover $f(P,x)$. Consider for
example the case when $P=[\0,\1]$ is a graded Eulerian poset of rank $n+1$. 
 Then, by \cite[(19)]{Stanley-GH},
\begin{equation}
\label{eq:[0,1]}
f([\0,\1],x)=x^{n+1}g([\0,\1),1/x)
\end{equation}
contains only terms of degree greater than $\lfloor n/2\rfloor$. Thus,
by Proposition~\ref{prop:lowerE}, we obtain $\st([\0,\1],x)=0$.
\end{remark}
Next we prove a generalization of the short toric variant of Fine's
formula~\eqref{eq:Fine} to all posets $P$ for which $\st(P,x)$ is
defined. Let $P$ be such a poset and let $n=\max\{\rank (p)\::\: p\in P\}$.
As for graded posets, for any $S\subseteq [1,n]$, we may define the number
$f_S=f_S(P)$ as the number of maximal chains in the $S$ rank-selected
subposet $P_S=\{u\in P\::\: \rank(u)\in S\}$. 

\begin{proposition}[Fine's formula] 
\label{prop:Fine-t}
\begin{equation}
\label{eq:Fine-t}
\st(P,x)=\sum_{S\subseteq [1,n]} f_S(P)\cdot \sum_{\lambda\in
\{-1,1\}^n\::\: S(\lambda)\supseteq S, n-2i_{\lambda}\geq 0}
(-1)^{n-i_{\lambda}+|S|}x^{n-2i_{\lambda}} 
\end{equation}
holds for all finite posets $P$ having a unique minimum element $\0$ and a
rank function $\rank: P\rightarrow {\mathbb N}$, satisfying $\rank(\0)=0$
and $n=\max\{\rank (p)\::\: p\in P\}$. 
\end{proposition} 
\begin{proof}
The statement may be shown in a very similar fashion to Fine's original
formula, the only difference being that the role of the intertwined
recurrence equations \eqref{eq:fgrec} and \eqref{eq:gfrec}
is taken over by the single recurrence given in Definition~\ref{def:stgen}, 
making the proof somewhat simpler. Thus we only outline the proof.
Introducing
$$
\st_f(S,n,x):=\sum_{\lambda\in
\{-1,1\}^n\::\: S(\lambda)\supseteq S, n-2i_{\lambda}\geq 0}
(-1)^{n-i_{\lambda}+|S|}x^{n-2i_{\lambda}},
$$
we need to show $\st(P,x)=\sum_{S\subseteq [1,n]} f_S(P)\cdot \st_f(S,n,x)$.
Observe first that the condition
$n-2i_{\lambda}\geq 0$ in the 
summation defining $\st_f(S,n,x)$ is equivalent to
$\lambda_1+\cdots+\lambda_n\geq 0$ and it may be dropped if we apply the
operator $U_{\geq 0}$ instead. In particular, we have
$$
\st_f(\emptyset,n,x)=U_{\geq 0}\left((x^{-1}-x)^n)\right),
$$
thus, by Definition~\ref{def:stgen}, we only need to prove
$$
U_{\geq 0}\left(\sum_{p\in P\setminus \{0\}} U_{\geq 1}\left(\st ([\0,p),x)
  (x-x^{-1})\right)(x^{-1}-x)^{n-\rank(p)}\right)=\sum_{\emptyset\neq 
S\subseteq [1,n]} f_S(P)\cdot \st_f(S,n,x).
$$
This statement may easily be shown by induction on $n$, using the fact
that, for $S\neq\emptyset$, we have
$$
\st_f(S,n,x)=U_{\geq 0}\left(U_{\geq 1}\left(\st_f(S\setminus\{\max
S\},\max S,x)(x-x^{-1})\right)(x^{-1}-x)^{n-\max S})\right)\quad\mbox{and}
$$ 
$$
f_S(P)=\sum_{p\in P: \rank(p)=\max S} f_{S\setminus \{\max S\}}([\0,p]). 
$$
\end{proof}
\begin{corollary}
\label{cor:stfulldeg}
The degree of $\st(P,x)$ equals $\max\{\rank(p)\::\: p\in P\}$ if and
only if 
$$\sum_{S\subseteq [1,n]} (-1)^{|S|} f_S(P)\neq 0.$$
\end{corollary}
It is worth noting that $\sum_{S\subseteq [1,n]} (-1)^{|S|} f_S(P)$ is
the reduced Euler characteristic of the {\em order complex} of
$P\setminus \{\0\}$: this is the simplicial complex, whose vertex set is
$P\setminus \{\0\}$ and whose faces are the chains in $P\setminus
\{\0\}$. It is not difficult to show by induction on $n$, using 
Definition~\ref{def:stgen}, that $\st(P,x)$ is additively symmetric for
every $P$. However, in light of Proposition~\ref{prop:lowerE} and
Remark~\ref{rem:lowerE}, it 
makes more sense to call the coefficients of the multiplicatively
symmetric variant of $\st(P,x)$ the ``toric $h$-vector of $P$'' when
$\st(P,x)$ has ``full degree''. According to
Corollary~\ref{cor:stfulldeg}, this is possible exactly when the 
reduced Euler characteristic of the order complex of
$P\setminus \{\0\}$ is not zero.

As it was the case with Fine's formula and Eulerian posets, the short
toric variant of Fine's formula above may be used to express $\st(P,x)$
in terms of the flag $h$-vector, whenever $\st(P,x)$ is defined.
For that purpose we may directly adapt the ideas present in the work of
Bayer and Ehrenborg~\cite[Section 7.4]{Bayer-Ehrenborg}. 
An outline may be found in the Appendix. 

\section{The short toric polynomial and the $cd$-index}
\label{sec:t-cd}

In this section we show how easy it is to compute the polynomial
$\st([\0,\1),x)$,  
associated to an Eulerian poset $[\0,\1]$, from its $cd$-index. 
We begin by expressing $\st([\0,\1),x)$ in terms of the
  $ce$-index. 

Using \eqref{eq:fL} and the binomial theorem
we may rewrite~\eqref{eq:Fine-t} as 
\begin{equation}
\label{eq:FineL-t}
\st([\0,\1),x)=\sum_{S\subseteq [1,n]} L_S \sum_{\lambda\in
    \{-1,1\}^n\::\: n-2i_{\lambda}\geq 0}
x^{n-2i_{\lambda}} (-1)^{n-i_{\lambda}+|S(\lambda)\setminus S|}. 
\end{equation}
Just like in~\cite[Section 7.4]{Bayer-Ehrenborg}, we represent 
each $\lambda\in\{-1,1\}^n$ by a lattice path $\Lambda(\lambda)$  
starting at $(0,0)$ and containing $(1,\lambda_i)$ as step $i$ for
$i=1,\ldots,n$. The resulting lattice path has $n$ steps, each step is a
northeast step $(1,1)$ or a southeast step $(1,-1)$. 
The condition $n-2i_{\lambda}\geq 0$ is equivalent to stating that we
only consider lattice paths whose right endpoint is on or above the
horizontal axis. As is usual in lattice-path combinatorics, we may 
use the ``reflection principle'' to match canceling terms, thus
simplifying the summation. For that purpose, let us introduce
$R(\lambda):=\{i\in [1,n]\::\: \lambda_1+\cdots+\lambda_i=0\}$. Note
that $R(\lambda)$ necessarily consists of even integers. We also say
that a set $S$ {\em   evenly contains} the set $R$ if $R\subseteq S$ and
$S\setminus R$ is the disjoint union of intervals of even cardinality.
\begin{proposition}
\label{prop:Finece-t}
Let $[\0,\1]$ be a graded Eulerian poset of rank $n+1$. Then we have 
$$
\st([\0,\1),x)=\sum_{S\subseteq [1,n]} L_S\cdot \st_{ce}(S,n,x).
$$
Here $\st_{ce} (S,n,x)$ is the total weight of all $\lambda\in\{-1,1\}^n$
such that $S$ evenly contains $R(\lambda)\cup (R(\lambda)-1)$ and
$\lambda_1+\cdots+\lambda_n\geq 0$ . The
weights are defined as follows:  
\begin{enumerate}
\item each $i\in [1,n]$ such that $\lambda_i=-1$ contributes a
  factor of $1/x$; 
\item each $i\in [1,n]$ such that $\lambda_i=1$ contributes a
  factor of $-x$; 
\item each $i\in S(\lambda)\setminus S$ contributes an additional factor of
  $-1$.
\end{enumerate}
\end{proposition}
\begin{proof}
If we remove the condition that $S$ evenly contains $R(\lambda)\cup
(R(\lambda)-1)$, we obtain an exact rephrasing of \eqref{eq:FineL-t}. We 
only need to show that the vectors $\lambda\in\{-1,1\}^n$ such that 
$S$ does not contain $R(\lambda)\cup (R(\lambda)-1)$ evenly, may be
matched into pairs whose weights cancel. 

Since $R(\lambda):=\{r_1,\ldots,r_k\}$ consists of even integers only,
the set $R(\lambda)\uplus (R(\lambda)-1)$ may be written as the disjoint union
\begin{equation}
\label{eq:R}
R(\lambda)\cup (R(\lambda)-1)=\biguplus_{i=1}^k [r_i-1,r_i].
\end{equation} 
We may assume that $S$ is an even set, otherwise $L_S=0$. Then 
$S:=\{s_1,s_2,\ldots,s_{2m}\}$ may be
written as 
\begin{equation}
\label{eq:S}
S=\biguplus_{j=1}^m [s_{2j-1},s_{2j-1}+1]. 
\end{equation} 
It is easy to see that $R(\lambda)\cup (R(\lambda)-1)$ is evenly
contained in $S$ if and only if each $[r_i-1,r_i]$ appearing on the
right hand side of \eqref{eq:R} equals some 
$[s_{2j-1},s_{2j-1}+1]$ on the right hand side of \eqref{eq:S}. 
Assume $S$ does not contain $R(\lambda)\cup (R(\lambda)-1)$ evenly and
let $i$ be the least index such that $[r_i-1,r_i]$ does not equal any 
$[s_{2j-1},s_{2j-1}+1]$. Then either $r_i-1\notin S$ or $r_i-1$ is the 
right end of some interval $[s_{2j-1},s_{2j-1}+1]$. Either way $|S\cap
[1,r_i-1]|$ is even. After setting $r_0:=0$, we may also state that 
$S\cap [1,r_{i-1}]$ has even cardinality: it is the empty set when
$i=1$, and it is the disjoint union of some intervals $[s_{2j-1},s_{2j-1}+1]$
when $i>1$. Combining the previous two observations, we obtain that 
$|S\cap [r_{i-1}+1,r_i-1]|$ is even. On the other hand,
$|[r_{i-1}+1, r_i-1]|$ is odd since $r_{i-1}$ and $r_i$ are even.
Therefore $|[r_{i-1}+1, r_i-1]\setminus S|$ is odd. Consider now the vector
$\tilde{\lambda}$ that corresponds to the lattice path obtained by
reflecting the part of the lattice path associated to $\lambda$ between
$(r_{i-1},0)$ and $(r_i,0)$. The reflected part has the same number of
northeast steps and southeast steps thus, up to sign, $\lambda$ and
$\tilde{\lambda}$ have the same weight. We also have
$R(\lambda)=R(\tilde{\lambda})$ and
$\tilde{\tilde{\lambda}}=\lambda$. Finally, because $|[r_{i-1}+1, 
  r_i-1]\setminus S|$ is odd, item (3) above implies that the weights of
$\lambda$ and $\tilde{\lambda}$ cancel.  
\end{proof}
We may use the ``reflection principle'' to compute the polynomials
$\st_{ce}(S,n,x)$ defined in Proposition~\ref{prop:Finece-t} 
even more efficiently, as the total weight of even fewer vectors
$\lambda \in \{1,1\}^n$.  
\begin{theorem}
\label{thm:Finece-t} 
Using a new weighting, for any $n\geq 1$ and any
$S\subseteq [1,n]$, the polynomial $\st_{ce}(S,n,x)$ also equals the
total weight of all $\lambda\in \{-1,1\}^n$ such that $S$ evenly
contains $R(\lambda)\cup (R(\lambda)-1)$ and $\lambda_1+\cdots
+\lambda_i\geq 0$ holds for all $i\in [1,n]$. The new weight of each
$\lambda$, satisfying the conditions stated above, is defined as follows:  
\begin{enumerate}
\item each $\lambda_i=-1$ contributes a factor of $-1/x$;
\item each $\lambda_i=1$ contributes a factor of $x$;
\item each element of $R(\lambda)$ contributes an additional factor of $2$. 
\end{enumerate}
\end{theorem}
\begin{proof}
Consider the way described in Proposition~\ref{prop:Finece-t}
to compute the polynomials $\st_{ce}(S,n,x)$ as the total weight of all 
$\lambda\in \{-1,1\}^n$. Let us say that $\nu\in\{-1,1\}^n$ is a 
reflection of $\lambda$ if there is an $r_i\in
R(\lambda)=\{r_1,\ldots,r_k\}$ such that the lattice path $\Lambda(\nu)$
corresponding to $\nu$ is obtained from the lattice path $\Lambda(\lambda)$
corresponding to $\lambda$ by reflecting the part of
$\Lambda(\lambda)$ between $(r_{i-1},0)$ and $(r_i,0)$. (As before,
we set $r_0:=0$.) An argument completely 
analogous to the proof of Proposition~\ref{prop:Finece-t} above shows
that $\lambda$ has the same weight as $\nu$ with respect to the weighting
defined there. Consider the
equivalence relation obtained by taking the transitive closure of 
the relation ``$\nu$ is a reflection of $\lambda$''. Clearly, the number
of elements in the equivalence class of $\lambda$ is $2^{|R(\lambda)|}$
and $\lambda$ being equivalent to $\nu$ implies $R(\lambda)=R(\nu)$ and 
$\lambda$ and $\nu$ have the same weight. Thus we may 
replace each $\lambda$ with the only vector $\lambda^+$ in its
equivalence class satisfying $\lambda^+_1+\cdots+\lambda^+_i\geq 0$ for
all $i\in [1,\max(R(\lambda^+))]$, at the expense of introducing an
additional factor of $2$, contributed by each $i\in R(\lambda^+)$. 
Obviously, each such class representative $\lambda^+$ exists uniquely. 
Observe next that each class representative $\lambda^+$ satisfies
$\lambda^+_1+\cdots+\lambda^+_i\geq 0$ for all $i$. This follows form 
$\lambda^+_1+\cdots+\lambda^+_n\geq 0$ and from the fact that the
associated lattice path of $\lambda^+$ can not cross the horizontal axis
after having reached the point $(\max(R(\lambda^+)),0)$. Thus we have
$S(\lambda^+)=[1,n]\setminus R(\lambda^+)$. Since $S$ contains 
$R(\lambda^+)$, the set $S(\lambda^+)=[1,n]\setminus R(\lambda^+)$ contains
$[1,n]\setminus S$ and $S(\lambda^+)\setminus S=[1,n]\setminus S$.
Thus, for each $\lambda^+$, item (3) in Proposition~\ref{prop:Finece-t} 
is equivalent to requiring that each $i\in [1,n]\setminus S$
contributes an additional factor of $-1$. Since $S$ is an even set,
$|[1,n]\setminus S|$ has the same parity as $n$. Therefore, for the 
class representatives $\lambda^+$, we may remove
item (3) in Proposition~\ref{prop:Finece-t} 
after changing the signs of the contributions in items (1) and (2).
\end{proof}
Theorem~\ref{thm:Finece-t} gains an even simpler form when we rephrase it
in terms of the $cd$-index.
\begin{theorem}
\label{thm:Finecd-t}
Let $[\0,\1]$ be a graded Eulerian poset of rank $n+1$. Then we have
$$
\st([\0,\1),x)=\sum_w [w] \Phi_{[\0,\1]}(c,d)\cdot \st(w,x).
$$
Here the summation runs over all $cd$-words $w$ of degree $n$. 
The polynomial $\st(w,x)$ is the total weight of all $\lambda\in 
\{-1,1\}^n$ such that the set of positions covered by letters $d$ equals  
$R(\lambda)\cup (R(\lambda)-1)$ and $\lambda_1+\cdots
+\lambda_i\geq 0$ holds for all $i\in [1,n]$. The weight of
each such $\lambda$ is defined as follows:  
\begin{enumerate}
\item each $\lambda_i=-1$ contributes a factor of $-1/x$;
\item each $\lambda_i=1$ contributes a factor of $x$;
\item each element of $R(\lambda)$ contributes an additional factor of $-1$. 
\end{enumerate}
\end{theorem}
\begin{proof}
Let us use the relation $d=\frac{1}{2}\cdot (c^2-e^2)$ to express each
$cd$-word $w$ as a linear combination 
\begin{equation}
\label{eq:cdce}
w=\sum_{w'} \alpha_{w,w'} w'
\end{equation}
of $ce$ words $w'$; and for each $ce$-word $w'$ let us denote by $S(w')$
the set of positions covered by a letter $e$ in $w'$. (The coefficient of
$w'$ in the $ce$-index is thus $L_{S(w')}$.) By
Proposition~\ref{prop:Finece-t} and by Theorem~\ref{thm:Finece-t}, it
suffices to show that the polynomial $\st(w,x)$ given in our statement
satisfies 
$$
\st(w,x)=\sum_{w'} \alpha_{w,w'} \st_{ce}(S(w'),n,x).
$$
Using the weighting introduced in Theorem~\ref{thm:Finece-t},
a vector $\lambda\in \{-1,1\}^n$ contributes 
to some $\st_{ce}(S(w'),n,x)$ only if the $S(w')$ evenly contains
$R(\lambda)\cup (R(\lambda)-1)$. Since the  factors $e^2$ only arise
from expanding $d$s, if $\lambda$ 
contributes to $\st(w,x)$ for some $cd$-word $w$, the set of positions
covered by $d$ must also evenly contain $R(\lambda)\cup
(R(\lambda)-1)$. We claim that, in order to have a nonzero contribution,
the two sets must be equal. Assume, by way of contradiction, that the
set of positions covered by $d$'s in $w$ properly contains
$R(\lambda)\cup (R(\lambda)-1)$. Due to the requirement of even
containment, this implies that there is a letter $d$ of $w$, covering the
positions $k$ and $k+1$, such that the set $\{k,k+1\}$ is disjoint from 
$R(\lambda)\cup (R(\lambda)-1)$. Consider the $ce$-words appearing in
the expansion \eqref{eq:cdce} of $w$. These all have the property that
the positions $k$ and $k+1$ are covered by either $c^2$ or $e^2$ and the
map $w'\mapsto \iota(w')$ replacing $c^2$ with $e^2$ and vice versa in
these positions is an involution. Our vector $\lambda$ contributes 
a nonzero weight to $\st_{ce}(S(w'),n,x)$ if and only if it contributes
a nonzero weight to $\st_{ce}(S(\iota(w')),n,x)$. Whenever this happens, the two
contributions cancel. 

We are left with considering vectors $\lambda\in\{-1,1\}^n$ such that
the set of positions covered by $d$'s in $w$ equals $R(\lambda)\cup
(R(\lambda)-1)$. For these, only the $ce$ word $w'$ obtained from
$w$ by replacing every $d$ with $-1/2\cdot e^2$ has the property that
the positions covered by $e$'s evenly contain $R(\lambda)\cup
(R(\lambda)-1)$. For this $w'$ we have 
$\alpha_{w,w'}=(-1/2)^{|R(\lambda)|}.$
The contribution of $\lambda$ to $\st(w,x)$ may be
described by the weighting above, considering that the factors 
of $2$, contributed by the elements of $R(\lambda)$ to $\st_{ce}(S(w'),n,x)$,
multiplied with $\alpha_{w,w'}$ leave us with a factor of $-1$ for
each element of $R(\lambda)$.
\end{proof}
Theorem~\ref{thm:Finecd-t} allows us to explicitly compute the
contribution $\st(w,x)$. Thus we obtain
the short toric equivalent of the result by Bayer and
Ehrenborg~\cite[Theorem 4.2]{Bayer-Ehrenborg}, expressing $f([\0,\1),x)$ in
terms of the $cd$-index.   
\begin{proposition}
\label{prop:Finecd-t}
The polynomial $\st(c^{k_1}dc^{k_2}\cdots c^{k_r}dc^{k},x)$ is zero if 
at least one of $k_1$, $k_2$, \ldots, $k_r$ is odd. If 
$k_1$, $k_2$, \ldots, $k_r$ are all even then 
$$
t(c^{k_1}dc^{k_2}\cdots c^{k_r}dc^{k},x)=
(-1)^{\frac{k_1+\cdots+k_r}{2}} C_{\frac{k_1}{2}} \cdots
C_{\frac{k_r}{2}} \QQ_{k}(x).
$$
Here $C_k=\frac{1}{k+1}\binom{2k}{k}$ is a Catalan number, and the
polynomials $\QQ_n(x)$ are given by $\QQ_0(x)=1$ and  
$$
\QQ_n(x):=\sum_{j=0}^{\left\lfloor \frac{n-1}{2}\right\rfloor} (-1)^j
\left(\binom{n-1}{j}-\binom{n-1}{j-1}\right) x^{n-2j}\quad\mbox{for
  $n\geq 1$.}
$$
\end{proposition}
\begin{proof}
To calculate $t(c^{k_1}dc^{k_2}\cdots c^{k_r}dc^{k_{r+1}},x)$ we must
sum the weight of all vectors $\lambda\in \{-1,1\}^n$ whose associated
lattice path is staying weakly above the horizontal axis, touching it
exactly at the points $(k_1+2,0)$, $(k_1+k_2+4,0)$, \ldots, $(k_1+\cdots
+k_r+2r,0)$. There is no such lattice path if at least one of $k_1$,
$k_2$, \ldots, $k_r$ is odd. Otherwise, we may select the lattice path
by independently selecting Dyck paths of length $k_1+2$, $k_2+2$,
\ldots, $k_r+2$ such that none of these Dyck paths touches the
horizontal axis between its start and end, and then we may independently
select a lattice path of length $k$ that stays strictly above the
horizontal axis. It is immediate from the definition of the weighting 
that a lattice path that ends at $(n,j)$ contributes a term $\pm x^j$. 
The total weight of the Dyck paths is 
$$
(-1)^{\frac{k_1+\cdots+k_r}{2}} C_{\frac{k_1}{2}} \cdots
C_{\frac{k_r}{2}},
$$
this needs to be multiplied by the total weight of all lattice paths of
length $k$ staying strictly above the horizontal axis. It is easy
to verify that this weight is $\QQ_{k}(x)$. 
\end{proof}
\begin{remark}
The polynomials $\QQ_n(x)$ are closely related to the polynomials
$Q_n(x)$ introduced by Bayer and Ehrenborg~\cite{Bayer-Ehrenborg}. 
They may be given by the formula
\begin{equation}
\QQ_n(x)=x^n Q_n(x^{-2}).
\end{equation}
\end{remark}
In analogy to the results thus far in this section, it is possible to
express $\st(P,x)$ of an arbitrary finite poset $P$ having a unique
minimum element $\0$ and a rank function $\rank: P\rightarrow {\mathbb
  N}$, in terms of its flag $h$-vector. This may be found in
the Appendix.

Besides yielding the short toric analogue of~\cite[Theorem
  4.2]{Bayer-Ehrenborg},  
Theorem~\ref{thm:Finecd-t} allows us to introduce two linear maps $\CC:
{\mathbb Q}[x]\rightarrow {\mathbb Q}[x]$ and $\DD: {\mathbb
  Q}[x]\rightarrow {\mathbb Q}[x]$ in such a way that, for any graded
Eulerian poset $[\0,\1]$, the polynomial $\st([\0,\1),x)$ may be
computed by substituting $\CC$ into $c$ and $\DD$ into $d$ in the
{\em reverse} of $\Phi_P(c,d)$ and applying the resulting linear operator 
to $1$. Note that the definitions and the result below are {\em independent of
the rank of $P$}.  
\begin{definition}
\label{def:CCDD}
We define $\CC: {\mathbb Q}[x]\rightarrow {\mathbb Q}[x]$ by
$$
\CC(x^n)=
\left\{
\begin{array}{ll}
x^{n+1}-x^{n-1} &\mbox{if $n\geq 2$,}\\
x^{n+1} &\mbox{if $n\in \{0,1\}$,}
\end{array}
\right.\quad
$$
and $\DD: {\mathbb Q}[x]\rightarrow {\mathbb Q}[x]$ by
$$
\DD(x^n)=\left\{
\begin{array}{ll}
-1 &\mbox{if $n=2$,}\\
1 &\mbox{if $n=0$,}\\
0 &\mbox{otherwise.}\\
\end{array}
\right.
$$
\end{definition}
\begin{theorem}
\label{thm:CCDD}
For any Eulerian poset $P=[\0,\1]$ we have 
$$
\st([\0,\1),x)=\Phi_{P}^{\rev}(\CC,\DD)(1)
$$
Here $\Phi_{P}^{\rev}(\CC,\DD)$ is obtained from $\Phi_P(c,d)$ by first taking
the reverse of each $cd$-monomial and then replacing each $c$ with $\CC$
and each $d$ with $\DD$.
\end{theorem}
\begin{proof}
By Theorem~\ref{thm:Finecd-t}, we only need to show that 
\begin{equation}
\label{eq:Finecd-t}
\st(c^{k_1}dc^{k_2}\cdots c^{k_r}dc^{k},x)=
\CC^k\DD \CC^{k_r}\DD\CC^{k_{r-1}}\DD\cdots \DD \CC^{k_1} (1)
\end{equation}
holds for any $cd$-word $w=c^{k_1}dc^{k_2}\cdots c^{k_r}dc^{k}$. This
may be shown by induction on the degree of $w$, the basis being 
$\st(\varepsilon, t)=1$ where $\varepsilon$ is the empty word. Assume now 
that \eqref{eq:Finecd-t} holds for some $cd$-word $w$ of degree $n$ above. By
Theorem~\ref{thm:Finecd-t}, each $\lambda\in\{-1,1\}^n$, contributing to
$\st(w,x)$ corresponds to a lattice path remaining weakly above the
horizontal axis such that each lattice point $(j,0)$ on the lattice path
satisfies that the positions $j$ and $j-1$ are covered by the same
letter $d$ in $w$. The weight contributed by $\lambda$ is $\pm
x^{\lambda_1+\cdots+\lambda_n}$ where $\lambda_1+\cdots+\lambda_n$ is
the height of the right end of the corresponding lattice path. To
compute $\st(wc,x)$, we must continue 
each lattice path contributing to $\st(w,x)$ with an additional step in
such a way that the resulting lattice path must end strictly above the
horizontal axis, and to multiply the weight of the lattice path with
the weight of the additional step. This amounts to applying 
$\CC$ to $\st(w,x)$, yielding 
$$\st(wc,x)=\CC\st(w,x).
$$
Similarly, to compute $\st(wd,x)$ we must continue each lattice
path contributing to $\st(w,x)$ with two steps in such a way that the
resulting lattice path ends on the horizontal axis. Clearly this can be
done to the lattice paths of positive length and of weight $\pm x^2$
or $1$ only. This amounts to applying $\DD$ to $\st(w,x)$, yielding 
$$
\st(wd,x)=\DD\st(w,x).
$$
\end{proof}

\section{Two useful bases}
\label{sec:bases}

Proposition~\ref{prop:Finecd-t} highlights the importance of the basis 
$\{\QQ_n(x)\}_{n\geq 0}$ of the vector space ${\mathbb Q}[x]$. In this
section we express the elements of the basis $\{x^n\}_{n\geq 0}$, as
well as the operators $\CC$ and $\DD$, in this new basis.  
We also find the analogous results for the basis $\{t_n(x)\}_{n\geq 0}$
where $t_n(x)$ is defined as $\st(B_{n+1},x)$ for the Boolean algebra
$\widehat{B}_{n+1}$ of rank $n+1$. This basis will be useful in proving the main
result of Section~\ref{sec:duals}, as well as in finding a very simple
formula connecting $\st(P,x)$ with $g(P,x)$. Both $\{\QQ_n(x)\}_{n\geq
  0}$ and $\{t_n(x)\}_{n\geq 0}$ are bases of ${\mathbb Q}[x]$ since
both sets contain exactly one polynomial of degree $n$ for each $n\geq
0$. 

\begin{proposition}
\label{prop:x-QQ}
For $n>0$ we have 
$$
x^n=\sum_{k=0}^{\left\lfloor \frac{n-1}{2}\right\rfloor}
\binom{n-1-k}{k} \QQ_{n-2k}(x).
$$
\end{proposition}
\begin{proof}
We claim that the right hand side is the total weight of all lattice paths
from $(0,0)$ to some point whose first coordinate is $n$, 
using northeast steps $(1,1)$ of weight $x$, southeast steps $(1,-1)$ of
weight $-1/x$, and ``long'' horizontal steps $(2,0)$ of weight $1$,
such that the entire path remains strictly above the horizontal axis.
Indeed, the 
first step of each such lattice path is a northeast step. After this step, we 
select the number $k\in [0,\lfloor (n-1)/2\rfloor]$ of
horizontal steps. By removing the horizontal steps from the list of
steps, we obtain a lattice path with altogether $(n-2k)$ nonhorizontal
steps that remain strictly above the horizontal axis.  
The total weight of all such lattice paths is $\QQ_{n-2k}(x)$.
For each such lattice path, there are
$\left(\binom{n-2k}{k}\right)=\binom{n-k-1}{k}$ ways to reinsert 
the $k$ horizontal steps after the first, second, \ldots, or $(n-2k)$-th
nonhorizontal step.

We only need to show that the total weight of the lattice paths
described above is $x^n$. The weight of the lattice path consisting of
$n$ northeast steps is $x^n$. We will show
that all other lattice paths may be arranged into pairs whose
weights cancel. For that purpose, consider any lattice path $\Gamma$
that contains at least one southeast step or at least one
horizontal step. If $\Gamma$ contains a southeast step before any
horizontal step, then this step must be immediately preceded by a 
northeast step. Let $i=i(\Gamma)$ be the least integer such that the
step starting at first coordinate $i$ is either a horizontal step,
or it is a northeast step immediately followed by a southeast step. Let
$\jmath (\Gamma)$ be the lattice path obtained from $\Gamma$ as follows. If the
step starting at first coordinate $i(\Gamma)$ is a horizontal
step, replace it with a northeast step followed by a southeast step. If
it is a northeast step followed by a southeast step, replace these two
steps with a single horizontal step. The bijection $\Gamma\mapsto
\jmath(\Gamma)$ is a fixed point free involution on the set of lattice
paths considered and the weight of $\Gamma$ is the negative of the
weight of $\jmath(\Gamma)$.  
\end{proof}

We may rewrite Proposition~\ref{prop:x-QQ} as 
\begin{equation}
\label{eq:x-QQ-e}
x^{2n}=\sum_{k=1}^{n} \binom{n-1+k}{n-k}\QQ_{2k} (x)
\end{equation}
for even powers of $x$ and as 
\begin{equation}
\label{eq:x-QQ-o}
x^{2n+1}=\sum_{k=0}^{n} \binom{n+k}{n-k}\QQ_{2k+1} (x)
\end{equation}
for odd powers of $x$. The coefficients appearing in equations 
\eqref{eq:x-QQ-e} and \eqref{eq:x-QQ-o} respectively are exactly the
coefficients of the {\em Morgan-Voyce polynomials} $B_n(x)$ and $b_n(x)$
respectively. They first appeared in the study of electrical
networks~\cite{Morgan-Voyce}, some of the other early references include
Swamy's work~\cite{Swamy1} and~\cite{Swamy2}. Another connection between
the toric $g$-polynomials of the cubes and the Morgan-Voyce
polynomials was noted by Hetyei~\cite{Hetyei-cubical}. 
\begin{corollary}
The linear transformation ${\mathbb Q}[x]\rightarrow {\mathbb Q}[x]$
given by $x^n\mapsto \QQ_n(x)$ takes $B_n(x^2)$ into $x^{2n}$ and 
$xb_n(x^2)$ into $x^{2n+1}$. 
\end{corollary}
Comparing Proposition~\ref{prop:Finecd-t} with \eqref{eq:Finecd-t}
yields the following consequence.
\begin{corollary}
\label{cor:CCDD-QQ}
The operators $\CC$ and $\DD$ are equivalently given by 
$$\CC(\QQ_n(x)):=\QQ_{n+1}(x)\quad \mbox{and}\quad 
\DD(\QQ_{n}(x))=
\left\{
\begin{array}{ll}
0 & \mbox {for odd $n$,}\\ 
(-1)^{n/2} C_{n/2}& \mbox {for even $n$.}\\ 
\end{array}
\right.
$$
\end{corollary} 
We now turn to the polynomials $t_n(x):=\st(\widehat{B}_{n+1},x)$. Stanley's
result~\cite[Proposition 2.1]{Stanley-GH} may be rewritten as
\begin{equation}
\label{eq:Boolean-t}
t_n(x)=\sum_{k=0}^{\left\lfloor\frac{n}{2}\right\rfloor} x^{n-2k}\quad
\mbox{for $n\geq 0$.}
\end{equation} 
Inverting the summation given in \eqref{eq:Boolean-t} yields
\begin{equation}
\label{eq:x-t}
x^n=
\left\{
\begin{array}{ll}
t_n(x)-t_{n-2}(x)&\mbox{if $n\geq 2$,}\\
t_n(x)&\mbox{if $n\in \{0,1\}$}.\\
\end{array}
\right.
\end{equation}
As an immediate consequence of Definition~\ref{def:CCDD} and
  \eqref{eq:x-t} we obtain 
\begin{equation}
\label{eq:CCDD-t}
\CC(t_n(x))=
t_{n+1}(x)-t_{n-1}(x) \quad\mbox{and}\quad
\DD(t_n(x))=\delta_{n,0} \quad\mbox{for $n\geq 0$}.
\end{equation}
Here we set $t_{-1}(x):=0$ and $\delta_{n,0}$ is the Kronecker delta
function. Finally, the most remarkable property of the basis
$\{t_n(x)\}_{n\geq 0}$ is its role in the following result connecting
the polynomials $g(P,x)$ and $\st(P,x)$.
\begin{proposition}
\label{prop:g-t}
Let $[\0,\1]$ be any Eulerian poset of rank $n+1$. Then 
\begin{equation}
\label{eq:st-c}
\st([\0,\1),x)=\sum_{k=0}^{\left\lfloor \frac{n}{2}\right\rfloor} c_k
  t_{n-2k}(x) 
\end{equation}
holds for some integers $c_0,c_1,\ldots, c_{\lfloor n/2\rfloor}$ if and
only if the same integers satisfy 
\begin{equation}
\label{eq:g-c}
g([\0,\1),x)=\sum_{k=0}^{\left\lfloor \frac{n}{2}\right\rfloor} c_k x^k.
\end{equation}. 
\end{proposition}
\begin{proof}
Assuming \eqref{eq:st-c}, equation \eqref{eq:g-c} is an
immediate consequence of Lemma~\ref{l:g-t} and the obvious identity 
\begin{equation}
U_{\geq 1}\left(\left(x-\frac{1}{x}\right)t_{m}(x)\right)=x^{m+1}.
\end{equation}
Conversely, assume \eqref{eq:g-c}. Then, it is easy to derive from
\eqref{eq:g-f} that 
$$
U_{\leq \lfloor n/2\rfloor}f([\0,\1),x)
=\sum_{k=0}^{\left\lfloor \frac{n}{2}\right\rfloor}
  \left(\sum_{i=0}^k c_i\right)  x^k,
$$
from which \eqref{eq:st-c} follows by the last part of \eqref{eq:t-f}.
\end{proof}

\section{The toric $h$-vector associated to an Eulerian dual simplicial poset}
\label{sec:duals}

In this section we use our formulas to compute 
the toric $h$-vector entries associated to and Eulerian {\em dual simplicial
poset} in terms of its $f$-vector. This question was mentioned by
Stanley~\cite{Stanley-GH}, citing an observation of Kalai. 

Given any graded poset $P$ of rank $n+1$, let us introduce $f_i$ for the
number of elements of rank $i+1$ in $P$ and call the resulting vector 
$(f_{-1},f_0,\ldots,f_n)$ the {\em $f$-vector} of $P$. If $P$ consists
of the faces of a polytope or simplicial complex, ordered by inclusion,
then $f_i$ is the number of faces of dimension $i$. A graded
poset $P$ is {\em simplicial} if for all $t\in P\setminus \{\1\}$, the
interval $[\0,t]$ is a Boolean algebra. For the toric $h$-polynomial
associated to a graded Eulerian simplicial poset $[\0,\1]$ of rank $n+1$
Stanley has the formula~\cite[Corollary 2.2]{Stanley-GH}    
\begin{equation}
\label{eq:th-s}
f([\0,\1),x)=\sum_{i=0}^n f_{i-1} (x-1)^{n-i}.
\end{equation}
In particular, if $[\0,\1)$ is also the face poset of a simplicial complex,
the toric $h$-vector of $P$ coincides with the $h$-vector of the
simplicial complex. 

A graded poset $P$ is {\em dual simplicial} if its dual $P^*$, obtained
by reversing the partial order, is a simplicial poset. Equivalently, for
all $t\in P\setminus \{\0\}$, the interval $[t,\1]$ is a Boolean
algebra. It was first observed by Kalai that the coefficients of the toric
$h$-polynomial of a dual simplicial
graded Eulerian poset $P$ depend only on the entries $f_i$ in the
$f$-vector $P$. This linear combination is not unique, and a simple
explicit formula was not known before. 
Given a dual simplicial Eulerian poset $P=[\0,\1]$, we will express the
toric $h$-polynomial coefficients of $[\0,\1)$ in terms of 
the toric $h$-vector $(h_0,\ldots,h_n)$ of the simplicial poset $(\0,\1]^*$.
Since $f_i(P)=f_{n-1-i}(P^*)$, by \eqref{eq:th-s} the entries of this
$h$-vector are given by 
$$
h_k=\sum_{i=k}^{n} (-1)^{i-k} \binom{i}{k} f_i(P).
$$
Thus, our formulas will show how the toric $h$-polynomial of a dual
simplicial poset $P$ depends on the numbers $f_i(P)$. It is worth noting
that in the special case when $P$ is the face lattice of a simple
polytope, the invariants $h_k$ given above are the entries of the
$h$-vector of the simple polytope.
We will apply the following result of Stanley~\cite[Theorem
  3.1]{Stanley-flag} to $P^*$.
\begin{theorem}[Stanley]
\label{thm:Stanley}
For each $n$, there exists $cd$-polynomials $\check{\Phi}^n_i$ (for
$i=0,\ldots,n-1$) such that the $cd$-index $\Phi_{P}(c,d)$ of any graded
Eulerian simplicial poset $P$ of rank $n+1$ may be written as 
$$
\Phi_{P}(c,d)=\sum_{i=0}^{n-1} h_i\cdot \check{\Phi}^n_i
$$ 
where $(h_0,\ldots,h_n)$ is the (toric) $h$-vector of $P$. 
\end{theorem}
The following description of the $cd$-polynomials $\check{\Phi}^n_i$ was
conjectured by Stanley~\cite[Conjecture 3.1]{Stanley-flag}
and shown by Hetyei~\cite[Theorem 2]{Hetyei-Andre}:   
\begin{theorem}
For each $n>0$ and each $0\leq i\leq n-1$, the polynomial
$\check{\Phi}^n_i$is the sum of the $cd$-variation monomials of all
augmented Andr\'e permutations $\pi$ of the set $[1,n+1]$,
ordered by the natural order, satisfying $\pi(n)=n-i$.
\end{theorem}
For our purposes, {\em augmented Andr\'e permutations} are most
conveniently defined as follows, see~\cite[Corollary 1]{Hetyei-Andre}.
\begin{definition}
\label{def:Andre}
Let $X$ be any linearly ordered set with $n$ elements. A {\em
  permutation} $\pi$ of $X$ is a word $\pi(1)\cdots 
\pi(n)$ in which each letter of $X$ occurs exactly once. An {\em
  augmented Andr\'e permutation} of $X$ is defined recursively as
follows:
\begin{itemize}
\item[(i)] For $n=0$ the empty word is an augmented
  Andr\'e permutation.
\item[(ii)] For $n=1$ the only permutation of $X$ is an augmented
  Andr\'e permutation.
\item[(iii)] For $n>1$, a permutation $\pi$ of $X$ is an augmented 
Andr\'e permutation if and only if for $m:=\pi^{-1}(\min X)$ the
permutations $\pi(1)\cdots \pi(m-1)$ and $\pi(m+1)\cdots \pi(n)$ are
augmented Andr\'e permutations and the letter $\max(X)$ belongs to
$\pi(m+1)\cdots \pi(n)$. Here $\min(X)$, respectively $\max(X)$, is the
least, respectively largest letter of $X$.
\end{itemize}
\end{definition}   
For further equivalent definitions and detailed bibliography on augmented
Andr\'e permutations we refer the reader to~\cite{Hetyei-Andre}. 
Given a permutation $\pi=\pi(1)\cdots\pi(n)$, the position
$i\in [1,n-1]$ is a {\em descent} if $\pi(i)>\pi(i+1)$,
otherwise it is an {\em ascent}. It is stated in one of the equivalent
definitions of augmented Andr\'e permutations~\cite[Definition 3]{Hetyei-Andre} 
  that every descent must be followed by an ascent. The {\em
    $cd$-variation monomial} of an augmented Andr\'e permutation 
$\pi(1)\cdots \pi(n)$ is obtained by placing a letter $d$ (of degree
$2$) to cover the positions $\{i,i+1\}$ for each descent $i$, and covering 
 all remaining positions in $[1,n-1]$ with a letter $c$ (of
 degree $1$). As an immediate consequence of Definition~\ref{def:Andre}
 we obtain that each augmented Andr\'e permutation must end with the
 largest letter. Furthermore, we have the following recurrence 
\begin{lemma}
\label{l:cphirec}
For $n\geq 2$, the noncommutative polynomials $\check{\Phi}^{n}_i$ satisfy 
$$
\check{\Phi}^{n}_i=
\left\{
\begin{array}{ll}
\Phi_{\widehat{B}_{n-1}}(c,d) d&\mbox{if $i=n-1$,}\\
c\check{\Phi}^{n-1}_i
+\sum\limits_{m=2}^{n-1}
\sum\limits_{j=0}^{\min(i,m-1)}\binom{i}{j}\binom{n-i-2}{m-1-j}
\Phi_{\widehat{B}_{m-1}}(c,d)d\check{\Phi}^{n-m}_{i-j}   
& \mbox{if $i<n-1$.}\\ 
\end{array}
\right.
$$
Here $\Phi_{\widehat{B}_{k}}(c,d)$ denotes the $cd$-index of the Boolean algebra of
rank $k$. 
\end{lemma}
\begin{proof}
An essentially equivalent statement may be found in
\cite[Proposition 7]{Hetyei-Andre} thus we only outline the proof of
this Lemma. An augmented Andr\'e permutation $\pi(1)\cdots \pi(n+1)$ of
$[1,n+1]$ must satisfy $\pi(n+1)=n+1$ and in
$\check{\Phi}^{n}_i$ we consider only Andr\'e permutations also
satisfying $\pi(n)=n-i$. Introducing
$m:=\pi^{-1}(1)$ we have $m=n$ if $i=n-1$ and $m\in [1,n-1]$
otherwise. The case $i=n-1$ follows immediately from the definitions. 
In the case $i<n-1$, term $c\check{\Phi}^n_i$ is contributed by the
augmented Andr\'e permutations $\pi(1)\cdots \pi(n+1)$ satisfying $m=1$,
the variable $j$ in the subsequent double sum stands for the size 
of the set $\{n-i+1,\ldots,n\}\cap \{\pi(1),\ldots,\pi(m-1)\}$. Besides
these $j$ elements, the set $\{\pi(1),\ldots,\pi(m-1)\}$ contains 
$m-1-j$ further elements of the set $\{2,\ldots,n-i-1\}$.
\end{proof}
As a consequence of Theorems~\ref{thm:CCDD} and Theorem~\ref{thm:Stanley} we
obtain the following statement. 
\begin{corollary}
\label{cor:duals-t}
The short toric polynomial $\st([\0,\1),x)$ associated to a graded dual
  simplicial Eulerian poset $P=[\0,\1]$ of rank $n+1$ may be written as
$$
\st([\0,\1),x)=\sum_{i=0}^{n-1} h_i \check{\Phi}^{n}_i(\CC,\DD)(1). 
$$ 
Here $(h_0,\cdots,h_n)$ is the toric $h$-vector of $(P\setminus\{\0\})^*$ and
$\check{\Phi}^{n}_i(\CC,\DD): {\mathbb Q}[x]\rightarrow
{\mathbb Q}[x]$ is the linear operator obtained by replacing each $c$
with $\CC$ and each $d$ with $\DD$ in the $cd$-polynomial
$\check{\Phi}^{n}_i$.  
\end{corollary}
Note that Theorem~\ref{thm:CCDD} involves reversing the order of the
letters when replacing $c$ with $\CC$ and $d$ with $\DD$. This change,
however, is offset by the fact that the $cd$-index of $P$ is obtained
from the $cd$-index of the simplicial poset $P^*$ by reversing the order
of the letters in each $cd$-monomial. We want to combine
Lemma~\ref{l:cphirec} with Corollary~\ref{cor:duals-t} to compute 
the short toric polynomial of a dual simplicial poset. This will lead to
a recurrence expressing the short toric polynomial contributions in terms of the
short toric polynomials of the Boolean algebras. 
As a special case of Corollary~\ref{cor:duals-t} we obtain the equation
\begin{equation}
\label{eq:B-d-t}
t_n(x)=\Phi_{\widehat{B}_{n+1}}(\CC,\DD)(1). 
\end{equation}
Here $\Phi_{\widehat{B}_{n+1}}(\CC,\DD)(1)$ is the operator obtained by replacing
each $c$ with $\CC$ and each $d$ with $\DD$ in the $cd$-index
$\Phi_{\widehat{B}_{n+1}}(c,d)$ of the Boolean algebra of rank $n+1$. 
\begin{theorem}
\label{thm:duals-t}
The short toric polynomial $\st([\0,\1),x)$ associated to a graded dual
  simplicial Eulerian poset $P=[\0,\1]$ of rank $n+1$ may be written as 
\begin{align*}
\st([\0,\1),x)&=h_0 (t_n(x)-(n-1)t_{n-2}(x))\\
&+\sum_{i=1}^{n-1} h_i \sum_{k=1}^{\left\lfloor
  \frac{n}{2}\right\rfloor}
\left(\binom{n-i}{k}\binom{i-1}{k-1}-\binom{n-i-1}{k}\binom{i}{k-1}\right)
t_{n-2k}(x).
\end{align*}
Here $(h_0,\cdots,h_n)$ is the toric $h$-vector of $(P\setminus\{\0\})^*$.
\end{theorem}
\begin{proof}
By Corollary~\ref{cor:duals-t} we only need to show that, for each 
$i\in \{0,\ldots,n-1\}$, the polynomial 
$t_{n,i}(x):=\check{\Phi}^{n}_i(\CC,\DD)(1)$ satisfies 
\begin{equation}
\label{eq:tni}
t_{n,i}(x)=
\left\{
\begin{array}{ll}
t_n(x)-(n-1)t_{n-2}(x) &\mbox{if $i=0$,}\\
\sum_{k=1}^{\left\lfloor
  \frac{n}{2}\right\rfloor}
\left(\binom{n-i}{k}\binom{i-1}{k-1}-\binom{n-i-1}{k}\binom{i}{k-1}\right)
t_{n-2k}(x) &\mbox{if $1\leq i\leq n-1$.}\\
\end{array}
\right.
\end{equation}
As before, in the formula above we use the convention $t_{-1}(x)=0$.
For small values of $n$ and $i$ the polynomials $t_{n,i}(x)$ are easy to
compute using \eqref{eq:CCDD-t} and the list of all polynomials
$\check{\Phi}^n_i$ satisfying $n\leq 5$, published by
Stanley~\cite{Stanley-flag}. The polynomials $t_{n,i}(x)$, up to $n=5$,
are listed in Table~\ref{tab:tni}. 
\begin{table}[h]
\begin{tabular}{l||l|l|l|l|l}
&$i=0$ & $i=1$ & $i=2$ & $i=3$ & $i=4$ \\
\hline
\hline
$n=1$ & $t_1(x)$&&&&\\
$n=2$ & $t_2(x)-t_0(x)$ & $t_0(x)$ &&&\\
$n=3$ & $t_3(x)-2t_1(x)$ & $t_1(x)$ &$t_1(x)$&&\\ 
$n=4$ & $t_4(x)-3t_2(x)$ & $t_2(x)-t_0(x)$ & $t_2(x)+t_0(x)$ &
$t_2(x)$&\\ 
$n=5$ & $t_5(x)-4t_3(x)$ & $t_3(x)-3t_1(x)$ & $t_3(x)+t_1(x)$ &
$t_3(x)+2t_1(x)$ & $t_3(x)$\\ 
\end{tabular}
\caption{The polynomials $t_{n,i}(x)$ for $n\leq 5$.}
\label{tab:tni}
\end{table}

We will prove \eqref{eq:tni} in several steps, by showing partial
statements using induction. In all such arguments our induction step
will depend on Lemma~\ref{l:cphirec}, combined with \eqref{eq:B-d-t}.
For $i=n-1$, the combination of these two statements yields 
\begin{equation}
\label{eq:tnn-1}
t_{n,n-1}(x)=t_{n-2}(x),
\end{equation}
whereas for $i<n-1$ we obtain
\begin{equation}
\label{eq:tnirec}
t_{n,i}(x)=\CC(t_{n-1,i}(x))+\sum\limits_{m=2}^{n-1}
\sum\limits_{j=0}^{\min(i,m-1)}
\binom{i}{j}\binom{n-i-2}{m-1-j}\DD(t_{n-m,i-j}(x))\cdot
t_{m-2}(x).
\end{equation}
Using these formulas we first prove the case $i=0$ in equation
\eqref{eq:tni}. The induction basis is $t_{1,0}(x)=t_1(x)$, see
Table~\ref{tab:tni}. For $n>1$ and the recurrence \eqref{eq:tnirec}
may be rewritten as 
$$
t_{n,0}(x)=\CC(t_{n-1,0}(x))+\sum\limits_{m=2}^{n-1}
\binom{n-2}{m-1}\DD(t_{n-m,0}(x))\cdot t_{m-2}(x).  
$$
This recurrence, combined with the induction hypothesis yields
\begin{align*}
t_{n,0}(x)=&\CC(t_{n-1}(x)-(n-2)t_{n-3}(x))\\ 
&+
\sum\limits_{m=2}^{n-1}
\binom{n-2}{m-1}\DD(t_{n-m}(x)-(n-m-1)\cdot t_{n-m-2}(x))\cdot t_{m-2}(x). 
\end{align*}
Using \eqref{eq:CCDD-t} this may be simplified to
$$
t_{n,0}(x)=t_n(x)-(n-1) t_{n-2}(x)+(n-2)
t_{n-4}(x)+\binom{n-2}{1}t_{n-4}(x)
=t_n(x)-(n-1) t_{n-2}(x),
$$
which is exactly what we wanted to prove.

Next we observe that all polynomials $t_{n,i}(x)$ are of the form 
$$
t_{n,i}(x)=\sum_{k=0}^{\left\lfloor\frac{n}{2}\right\rfloor}
\tau_{n,i,k}\cdot t_{n-2k}(x)
$$
for some integers $\tau_{n,i,k}\in {\mathbb Z}$. This may be easily
shown by induction on $n$, using the formulas \eqref{eq:tnn-1} and
\eqref{eq:tnirec} and the following two observations:
\begin{enumerate} 
\item When we apply the operator $\CC$ to an integer linear combination
  of polynomials $t_l(x)$ whose indices are all of the same parity, we
  obtain an integer linear combination of polynomials $t_l(x)$, and all
  indices are of the opposite parity.
\item To obtain a nonzero contribution, the input $t_{n-m,i-j}(x)$ of
  the operator $\DD$  must contain $t_{0}(x)$; thus by the induction hypothesis 
$m$ must have the same parity as $n$. The operator $\DD$ sends integers
  into integers.
\end{enumerate}  
Next we rewrite the formulas \eqref{eq:tnn-1} and
\eqref{eq:tnirec} as recurrences for the coefficients $\tau_{n,i,k}$. 
To facilitate this task, we extend the definition of 
$\tau_{n,i,k}$ to $n=0$ and to $i\geq n$ (where $n\geq 0$) by setting
$\tau_{n,i,k}=0$ if $n=0$ or $i\geq n$. For $n\geq 2$ and $i=n-1$,
\eqref{eq:tnn-1} may be rewritten as 
\begin{equation}
\label{eq:taunn-1}
\tau_{n,n-1,k}=\delta_{k,1},
\end{equation}
where $\delta_{k,1}$ is the Kronecker delta function. When rewriting 
\eqref{eq:tnirec} recall that, as observed above,
$\DD(t_{n-m,i-j}(x))=0$ unless $m$ is of the form $n-2s$ for some $s$
satisfying $1\leq s\leq \lfloor (n-2)/2\rfloor$. Thus \eqref{eq:tnirec} may be
rewritten as  
$$
t_{n,i}(x)=\CC(t_{n-1,i}(x))
+\sum\limits_{s=1}^{\left\lfloor\frac{n-2}{2}\right\rfloor} 
\sum\limits_{j=0}^{\min(i,n-2s-1)}
\binom{i}{j}\binom{n-i-2}{n-2s-1-j}\DD(t_{2s,i-j}(x))\cdot
t_{n-2s-2}(x).
$$
Since $n-2s-2=n-2k$ is equivalent to $s=k-1$, comparing coefficients of
$t_{n-2k}$ on both sides of the equation above yields
\begin{equation}
\label{eq:taunirec}
\tau_{n,i,k}=\tau_{n-1,i,k}-\tau_{n-1,i,k-1}+
\sum\limits_{j=0}^{\min(i,n-2k+1)}
\binom{i}{j}\binom{n-i-2}{n-2k+1-j}\tau_{2k-2,i-j,k-1}\quad\mbox{for $i<n-1$.}
\end{equation}
$\tau_{2k-2,i-j,k-1}=0$ when $k-1=2k-2=0$ or $i-j\geq 2k-2$.
We only need to show the case $i>0$ of \eqref{eq:tni}, which is
equivalent to
\begin{equation}
\label{eq:tauni}
\tau_{n,i,k}=\binom{n-i}{k}\binom{i-1}{k-1}-\binom{n-i-1}{k}\binom{i}{k-1}
\quad\mbox{for $i>0$.}
\end{equation}
We prove this statement by induction on $n$, using the polynomials
listed in Table~\ref{tab:tni} as our induction basis. Observe that,
for $i=n-1$, \eqref{eq:taunn-1} gives the same $\delta_{k,1}$ as
\eqref{eq:tauni}. Thus we only need to show the validity of
\eqref{eq:tauni} in the case when $n\geq 2$ and $i<n-1$, assuming the
validity of \eqref{eq:tni} and, in particular, \eqref{eq:tauni} for all smaller
values of $n$. Since $k\leq \lfloor n/2\rfloor$, we have $2k-2<n$, thus
we may use \eqref{eq:taunirec} to compute $\tau_{n,i,k}$. 

{\bf\noindent Case 1.} $k=0$. 
In this case \eqref{eq:taunirec} may be simplified to 
$\tau_{n,i,0}=\tau_{n-1,i,0}$. Repeated application of this recurrence
yields $\tau_{n,i,0}=\tau_{i+1,i,0}$ which equals $0$ by
\eqref{eq:taunn-1}. Formula \eqref{eq:tauni} also gives $0$.  

{\bf\noindent Case 2.} $k=1$. 
In this case \eqref{eq:taunirec} may be simplified to 
$\tau_{n,i,1}=\tau_{n-1,i,1}-\tau_{n-1,i,0}$. By the already shown
previous case (and also by our induction hypothesis) we have
$\tau_{n-1,i,0}=0$, thus we obtain the recurrence  
$\tau_{n,i,1}=\tau_{n-1,i,1}$. Repeated application of this recurrence
yields $\tau_{n,i,1}=\tau_{i+1,i,1}$ which equals $1$ by
\eqref{eq:taunn-1}. Formula \eqref{eq:tauni} also gives $1$.  

{\bf\noindent Case 3.} $k=2$. 
By Table~\ref{tab:tni} we have $\tau_{2,0,1}=-1$ and $\tau_{2,1,1}=1$.
Thus \eqref{eq:taunirec} may be simplified to
$$
\tau_{n,i,2}=\tau_{n-1,i,2}-\tau_{n-1,i,1}+
\binom{i}{i-1}\binom{n-i-2}{n-i-2}-
\binom{i}{i}\binom{n-i-2}{n-3-i}
\quad\mbox{for
  $i<n-1$.}
$$
This recurrence equals to the recurrence obtained by substituting $k=2$
into \eqref{eq:taunirecsimp} in Case 4 below. The rest of the proof
of this case is identical to the proof of Case 4.

{\bf\noindent Case 4.} $k\geq 3$. Consider the term
$\tau_{2k-2,i-j,k-1}$ in the sum on the right hand side of
\eqref{eq:taunirec}. If $i-j=0$ then, by the already shown first part of
\eqref{eq:tni} we have
$$
t_{2k-2,0}(x)=t_{2k-2}(x)-(2k-3) t_{2k-4}(x),
$$ 
implying $\tau_{2k-2,0,k-1}=0$, since $k-1\geq 2$. If $i-j>0$ then, by
our induction hypothesis we have 
$$
\tau_{2k-2,i-j,k-1}=
\binom{2k-2-(i-j)}{k-1}\binom{i-j-1}{k-2}
-\binom{2k-3-(i-j)}{k-1}\binom{i-j}{k-2}. 
$$
Here we have $\binom{2k-2-(i-j)}{k-1}=\binom{2k-3-(i-j)}{k-1}=0$ unless 
$i-j\leq k-1$ and $\binom{i-j-1}{k-2}=\binom{i-j}{k-2}=0$ unless
$i-j\geq k-2$. As a consequence, $\tau_{2k-2,i-j,k-1}=0$ unless $i-j\in
\{k-2,k-1\}$. Direct substitution into the above formula
shows $\tau_{2k-2,k-2,k-1}=-1$ and
$\tau_{2k-2,k-1,k-1}=1$. Therefore, the recurrence \eqref{eq:taunirec}
may be simplified to 
\begin{equation}
\label{eq:taunirecsimp}
\tau_{n,i,k}=\tau_{n-1,i,k}-\tau_{n-1,i,k-1}+
\binom{i}{i-k+1}\binom{n-i-2}{n-k-i}
-
\binom{i}{i-k+2}\binom{n-i-2}{n-k-1-i},
\end{equation} 
Substituting our induction hypothesis for $\tau_{n-1,i,k-1}$ and 
using the symmetry of the binomial coefficients yields
\begin{align*}
\tau_{n,i,k}=&
\tau_{n-1,i,k}-
\binom{n-i-1}{k-1}\binom{i-1}{k-2}
+\binom{n-i-2}{k-1}\binom{i}{k-2}\\
&+\binom{i}{k-1}\binom{n-i-2}{k-2}
-\binom{i}{k-2}\binom{n-i-2}{k-1}\\
=&\tau_{n-1,i,k}-
\binom{n-i-1}{k-1}\binom{i-1}{k-2}
+\binom{i}{k-1}\binom{n-i-2}{k-2}.
\end{align*}
Substituting our induction hypothesis for $\tau_{n-1,i,k}$ yields
\begin{align*}
\tau_{n,i,k}=& \binom{n-i-1}{k}\binom{i-1}{k-1}-
\binom{n-i-2}{k}\binom{i}{k-1}
+\binom{n-i-2}{k-2}\binom{i}{k-1}\\
&-\binom{n-i-1}{k-1}\binom{i-1}{k-2}.\\
\end{align*}
By Pascal's identity, we may replace $\binom{i-1}{k-2}$ with
$\binom{i}{k-1}-\binom{i-1}{k-1}$ and $\binom{n-i-2}{k-2}$ with
$\binom{n-i-1}{k-1}-\binom{n-i-2}{k-1}$  in the last equation. Thus we obtain
\begin{align*}
\tau_{n,i,k}=& \binom{n-i-1}{k}\binom{i-1}{k-1}-
\binom{n-i-2}{k}\binom{i}{k-1}
+\binom{n-i-1}{k-1}\binom{i}{k-1}\\
&-\binom{n-i-2}{k-1}\binom{i}{k-1}-\binom{n-i-1}{k-1}\binom{i}{k-1}+\binom{n-i-1}{k-1}\binom{i-1}{k-1}.\\
\end{align*}
After collecting terms by factors of $\binom{i-1}{k-1}$ and 
$\binom{i}{k-1}$, respectively, and using Pascal's identity,
we obtain \eqref{eq:tauni}.
\end{proof}
An important equivalent form of Theorem~\ref{thm:duals-t} is the
following statement.
\begin{proposition}
\label{prop:duals-tg}
Let $[\0,\1]$ be a graded dual simplicial Eulerian poset of rank $n+1$
and let $(h_0,\cdots,h_n)$ be the toric $h$-vector of $(\0,\1]^*$. 
Then we have
$$
\st([\0,\1),x)=h_0 t_n(x) +\sum_{i=1}^{\left\lfloor \frac{n}{2}\right\rfloor}
  (h_i-h_{i-1})
\sum_{k=1}^{\min \{i,n-i\}} \frac{n+1-2i}{k}\binom{n-i}{k-1}\binom{i-1}{k-1}
t_{n-2k}(x). 
$$
\end{proposition}
\begin{proof}
As an immediate consequence of Theorem~\ref{thm:duals-t}, we have 
$[t_{n}(x)] \st([\0,\1),x)=h_0$. We only need to prove that the
  coefficient of $t_{n-2k}(x)$ in $\st([\0,\1),x)$ is correctly stated
    for $k\geq 1$. 
In terms of the coefficients $\tau_{n,i,k}$,
    introduced in the proof of Theorem~\ref{thm:duals-t}, we have
$$
[t_{n-2k}(x)] \st([\0,\1),x)=\sum_{i=0}^{n-1} h_i\cdot \tau_{n,i,k}.
$$
Using the fact that $(h_0,\ldots,h_n)$ is the toric $h$-vector of the dual
of $[\0,\1]$ and thus satisfies $h_i=h_{n-i}$ for $i=1,2,\ldots,n,$ we
may rewrite the above equation as 
$$
[t_{n-2k}(x)] \st([\0,\1),x)=h_0\cdot \tau_{n,0,k}+\sum_{i=1}^{\lfloor n/2 \rfloor}
  h_i\cdot (\tau_{n,i,k}+\tau_{n,n-i,k})
-\delta_{\lfloor n/2\rfloor,\lceil n/2 \rceil}
\cdot h_{\lfloor n/2\rfloor}\cdot \tau_{n,\lfloor n/2\rfloor,k}. 
$$
Here $\delta_{\lfloor n/2\rfloor,\lceil n/2 \rceil}$ is the Kronecker
delta, and adding the last term represents subtracting $h_{n/2}\cdot
\tau_{n,n/2,k}$ exactly when $n$ is even. Rewriting the right hand side as a
linear combination of $h_0$, $h_1-h_0$, \ldots, $h_{\lfloor
  n/2\rfloor}-h_{\lfloor n/2\rfloor-1}$ yields
$$
[t_{n-2k}(x)] \st([\0,\1),x)=h_0\cdot \left(\tau_{n,0,k}
+\sum_{j=1}^{n-1} \tau_{n,j,k}\right)
+\sum_{i=1}^{\lfloor n/2 \rfloor}
  (h_i-h_{i-1})\cdot \sum_{j=i}^{n-i} \tau_{n,j,k}. 
$$
The statement now follows from \eqref{eq:tni} and from the fact that 
\begin{align*}
\sum_{j=i}^{n-i} \tau_{n,j,k}
&=
\sum_{j=i}^{n-i}
\left(\binom{n-j}{k}\binom{j-1}{k-1}-\binom{n-j-1}{k}\binom{j}{k-1}\right)\\ 
&=\sum_{j=i-1}^{n-i-1} \binom{n-j-1}{k}\binom{j}{k-1}
-\sum_{j=i}^{n-i} \binom{n-j-1}{k}\binom{j}{k-1}\\
&=\binom{n-i}{k}\binom{i-1}{k-1}
-\binom{i-1}{k}\binom{n-i}{k-1}
=\frac{n+1-2i}{k} \binom{n-i}{k-1}\binom{i-1}{k-1}
\end{align*}
holds for $i=1,2,\ldots, \lfloor n/2\rfloor$. Note that, applying the
last equation to $i=1$ implies
$$
\tau_{n,0,k}+\sum_{j=1}^{n-1} \tau_{n,j,k}=\tau_{n,0,k}+\frac{n-1}{k}
\binom{n-1}{k-1}\binom{0}{k-1}=
\left\{ 
\begin{array}{ll}
-(n-1)+(n-1)&\mbox{for $k=1$,}\\
0+0 &\mbox{for $k\geq 2$.}
\end{array}
\right.
$$
In either case, $h_0$ contributes zero to $t_{n-2k}(x)$ for $k\geq 1$.
\end{proof}

Using Proposition~\ref{prop:g-t}, Theorem~\ref{thm:duals-t} and
Proposition~\ref{prop:duals-tg}, respectively, may be rewritten as the
formulas stated in the next two corollaries.
\begin{corollary}
\label{cor:duals-g}
Let $[\0,\1]$ be a graded dual simplicial Eulerian poset of rank $n+1$
and let $(h_0,\cdots,h_n)$ be the toric $h$-vector of $(\0,\1]^*$. 
Then 
\begin{align*}
g([\0,\1),x)&=h_0 (1-(n-1)x)\\
&+\sum_{i=1}^{n-1} h_i \sum_{k=1}^{\left\lfloor
  \frac{n}{2}\right\rfloor}
\left(\binom{n-i}{k}\binom{i-1}{k-1}-\binom{n-i-1}{k}\binom{i}{k-1}\right)
x^k.
\end{align*}
\end{corollary}
\begin{corollary}
\label{cor:duals-gg}
Let $[\0,\1]$ be a graded dual simplicial Eulerian poset of rank $n+1$
and let $(h_0,\cdots,h_n)$ be the toric $h$-vector of $(\0,\1]^*$. 
Then 
$$
g([\0,\1),x)=h_0+\sum_{i=1}^{\left\lfloor \frac{n}{2}\right\rfloor}
  (h_i-h_{i-1}) 
\sum_{k=1}^{\min \{i,n-i\}} \frac{n+1-2i}{k}\binom{n-i}{k-1}\binom{i-1}{k-1}
x^k. 
$$
\end{corollary}
The most important consequence of Corollary~\ref{cor:duals-gg} is the
following.
\begin{corollary} 
\label{cor:glbcs}
Let $[\0,\1]$ be a graded dual simplicial Eulerian poset of rank $n+1$.
If the toric $h$-vector $(h_0,\ldots,h_n)$ of $(\0,\1]^*$ satisfies
$h_0\leq h_1\leq \cdots \leq h_{\lfloor n/2\rfloor}$, then
$f([\0,\1],x)$ has nonnegative coefficients.
\end{corollary}
Indeed, by Corollary~\ref{cor:duals-gg} above, $g([\0,\1),x)$ has
  nonnegative coefficients and the statement follows from
  \eqref{eq:[0,1]}. 
\begin{example}
Let $[\0,\1]$ be the face lattice of an $n$-dimensional simple
polytope ${\mathcal P}$. This is an Eulerian dual simplicial poset of
rank $n+1$; its 
dual is the face lattice of a simplicial polytope. The coefficients 
of $f([\0,\1],x)$ form the toric $h$-vector of ${\mathcal
  P}$. By Corollary~\ref{cor:glbcs}, the fact that the toric
$h$-vector of ${\mathcal P}$ has nonnegative entries is a
consequence of the Generalized Lower Bound Theorem~\cite{Stanley-NF} for
simplicial polytopes. 
\end{example}
\begin{remark}
As pointed out by Christian Krattenthaler~\cite{Kratt-pers}, the
coefficient 
$$
\frac{n+1-2i}{k}\binom{n-i}{k-1}\binom{i-1}{k-1}
=\binom{n+1-i}{k}\binom{i-1}{k-1}-\binom{n-i}{k-1}\binom{i}{k} 
$$
of $x^k$ in the contribution of $h_{i}-h_{i-1}$ to the right hand side
of Corollary~\ref{cor:duals-gg}
is the number of lattice paths from $(0,0)$ to $(n+1-i,i)$ using $N$ steps
$(1,0)$ and $E$ steps $(0,1)$ with exactly $k$ NE turns, such that the
lattice path does not cross the line $y=x$ and has its last NE turn on
the line $y=i$. In particular, in the case when $n=2i$, the coefficients
$N(i,k)=\binom{i-1}{k-1}\binom{i}{k-1}/{k}$, 
contributed by $h_{\lfloor n/2\rfloor}-h_{\lfloor n/2\rfloor-1}$ are
the {\em Narayana numbers}. The above description may be shown
with a slight modification of the proof of~\cite[Theorem
  3.4.1]{Kratt}. A direct application of~\cite[Theorem 3.4.1]{Kratt}
yields that, for $i\geq (n-1)/2$, the coefficients   
$$
\binom{n-i}{k}\binom{i-1}{k-1}-\binom{n-i-1}{k}\binom{i}{k-1}
=
\binom{n-i-1}{k-1}\binom{i}{k-1}-\binom{n-i}{k}\binom{i-1}{k-2},
$$
contributed by $h_i$ in Corollary~\ref{cor:duals-g},
count all lattice paths from $(0,0)$ to $(i,n-i-1)$ using $N$ steps
and $E$ steps with exactly $k-1$ NE turns, such that the
lattice path does not cross the line $y=x$. In particular, the
coefficients of the contributions of $h_{\lfloor n/2\rfloor}$ and
$h_{\lceil n/2\rceil}$ are again the Narayana numbers. For a detailed
bibliography of the Narayana numbers we refer to sequence A001263 in the
On-Line Encyclopedia of Integer Sequences~\cite{Sloane}. For $i<(n-1)/2$,
the contributions of $h_i$ in Corollary~\ref{cor:duals-g} also contain
negative coefficients, it seems hard to see a pattern of signs.
\end{remark}
Motivated by the well-known example of the cube (see
Example~\ref{ex:cube} below), we rewrite Corollary~\ref{cor:duals-g} in
the basis $\{(x-1)^k\}_{k\geq 0}$.
\begin{proposition}
\label{prop:duals-g}
Let $[\0,\1]$ be a graded dual simplicial Eulerian poset of rank $n+1$,
and let $(h_0,\cdots,h_n)$ be the toric $h$-vector of $(\0,\1]^*$.
Then we have
\begin{align*}
g([\0,\1),x)&=h_0 (-n+2-(n-1)(x-1))\\
&+\sum_{i=1}^{n-1} h_i \sum_{k=0}^{\left\lfloor
  \frac{n}{2}\right\rfloor}
\left(\binom{n-i}{k}\binom{n-k-1}{i-k}
-\binom{n-i-1}{k}\binom{n-k-1}{i+1-k}\right) (x-1)^k.
\end{align*}
\end{proposition}
\begin{proof}
The contribution of $h_0$ is the same in the above formula and in 
Corollary~\ref{cor:duals-g}. We only need to verify that each 
$h_i$ has the same contribution when $i>0$. By the binomial
theorem and the symmetry of the binomial coefficients,
 when we expand the contribution of $h_i$ given in
Corollary~\ref{cor:duals-g} in the basis $\{(x-1)^k\}_{k\geq 0}$, the 
coefficient of $(x-1)^k$ is 
\begin{equation}
\label{eq:sigmadef}
\sigma_{n,i,k}:=
\sum_{j=k}^{\left\lfloor\frac{n}{2}\right\rfloor} \binom{j}{k}
\left(\binom{n-i}{j}\binom{i-1}{i-j}
-\binom{n-i-1}{j}\binom{i}{i+1-j}\right).
\end{equation}
Using the identity
$\binom{n}{m}\binom{m}{k}=\binom{n}{k}\binom{n-k}{m-k}$ we obtain
\begin{align*}
\sigma_{n,i,k}=&
\binom{n-i}{k}\sum_{j=k}^{\left\lfloor\frac{n}{2}\right\rfloor} 
\binom{n-i-k}{j-k}\binom{i-1}{i-j}\\
&-\binom{n-i-1}{k}\sum_{j=k}^{\left\lfloor\frac{n}{2}\right\rfloor} 
\binom{n-i-1-k}{j-k}\binom{i}{i+1-j}.\\
\end{align*}
The statement now follows from the Chu-Vandermonde identity. Note that
the terms in the last two sums are zero unless $j\leq \min(i+1,n-i-1)$ which
is a more stringent condition than $j\leq \lfloor n/2\rfloor$.
\end{proof}
Using Pascal's identity and the symmetry of
the binomial coefficients, Proposition~\ref{prop:duals-g} yields
\begin{align*}
\sigma_{n,i,k}=&
\binom{n-i-1}{k}\binom{n-k-1}{n-i-1}
+\binom{n-i-1}{k-1}\binom{n-k-1}{n-i-1}\\
&-\binom{n-i-2}{k}\binom{n-k-1}{n-i-2}
-\binom{n-i-2}{k-1}\binom{n-k-1}{n-i-2}.\\
\end{align*}
for $i>0$ and the coefficients $\sigma_{n,i,k}$ defined in
\eqref{eq:sigmadef}. This may be rewritten as
\begin{equation}
\label{eq:sigma}
\sigma_{n,i,k}=
\binom{n-k-1}{k}p(n-2k-1,n-i-k-1)
+\binom{n-k-1}{k-1}p(n-2k,n-i-k)\quad\mbox{for $i>0$,}
\end{equation}
where we follow the notation of Bayer and
Ehrenborg~\cite{Bayer-Ehrenborg} by setting
$p(n,k):=\binom{n}{k}-\binom{n}{k-1}$.  
\begin{example}
\label{ex:cube}
Consider the face lattice $\widehat{L_n}$ of an $n$-dimensional cube. This is a
graded, dual simplicial Eulerian poset of rank $n+1$, satisfying
$h_i=\binom{n}{i}$. Thus  \eqref{eq:sigma} and the Chu-Vandermonde
identity yields
$$
\sum_{i=1}^{n-1} \binom{n}{i} \sigma_{n,i,k}=\binom{n-k-1}{k}p(2n-2k-1,n-k-1)
+\binom{n-k-1}{k-1}p(2n-2k,n-k)
$$
for $k\geq 2$. Note that the assumption $k\geq 2$ is necessary to make
sure that the omitted substitution $i=0$ yields only zero terms in 
all applications of the Chu-Vandermonde identity. Using the fact that 
$p(2n-2k-1,n-k-1)$ and $p(2n-2k,n-k)$ both equal the Catalan number
$C_{n-k}$ we obtain 
$$
\sum_{i=1}^{n-1} \binom{n}{i} \sigma_{n,i,k}=\binom{n-k}{k}C_{n-k}.
$$
Similarly, for $k=0$ and $k=1$ respectively, we obtain
$$ \sum_{i=1}^{n-1} \binom{n}{i}
\sigma_{n,i,0}=\binom{n}{0}C_{n}+n-2 \quad\mbox{and}\quad
\sum_{i=1}^{n-1} \binom{n}{i}
\sigma_{n,i,1}=\binom{n-1}{1}C_{n-1}+n-1.
$$
The additional terms on the right hand sides account for the omitted
substitutions $i=0$ on the left hand sides. Considering that $h_0=1$
contributes $-n+2-(n-1)(x-1)$, we find
\begin{equation}
\label{E_gL1}
g(L_n,x)=\sum_{k=0}^{\lfloor n/2\rfloor} 
\binom{n-k}{k}C_{n-k}(x-1)^k. 
\end{equation}
Here $L_n$ is obtained from $\widehat{L_n}$ by removing its maximum
element. It was noted in~\cite[Lemma 3.3]{Hetyei-cubical} that the above
formula  is equivalent to Gessel's result~\cite[Proposition
  2.6]{Stanley-GH}, stating 
\begin{equation}
\label{E_gL}
g(L_n,x)=\sum_{k=0}^{\lfloor n/2\rfloor} \frac{1}{n-k+1}
\binom{n}{k}\binom{2n-2k}{n} (x-1)^k. 
\end{equation}
The first combinatorial interpretation of the right hand side of
\eqref{E_gL} is due to Shapiro~\cite[Ex. 3.71g]{Stanley-EC1} the proof
of which was published by Chan~\cite[Proposition 2]{Chan}.  
\end{example}

\setcounter{section}{1} 
\renewcommand{\thesection}{\Alph{section}}
\renewcommand{\thesubsection}{\thesection.\arabic{subsection}}
\section*{Appendix}

Here we outline how to derive a formula expressing $\st(P,x)$
in terms of the flag $h$-vector of an arbitrary poset $P$ that has a
unique minimum element $\0$ and a rank function $\rank: P\rightarrow
{\mathbb N}$, satisfying $\rank(\0)=0$ and $\max\{\rank(p)\::\: p\in P\}=n$. 
There is no need nor place to use the reflection principle this
time, the line of thought present here is essentially an
adaptation of the work of Bayer and Ehrenborg~\cite[Section
  7.4]{Bayer-Ehrenborg}. 

In analogy to the proof of~\cite[Theorem 7.14]{Bayer-Ehrenborg} we
begin with substituting \eqref{eq:h-f} into \eqref{eq:Fine-t} and
rearranging. This yields  
$$
\st(P,x)=\sum_{\lambda\in \{-1,1\}^n \::\:
  n-2i_{\lambda}\geq 0} x^{n-2 i_{\lambda}}
(-1)^{|S(\lambda)|+n-i_{\lambda}} h_{S(\lambda)}.
$$
This equation may be rewritten as
\begin{equation}
\label{eq:st-h}
\st(P,x)=\sum_{S\subseteq [1,n]} h_S\cdot \st_h(S,x),
\end{equation}
where $\st_h(S,x)$ is the total weight of 
all $\lambda\in\{-1,1\}^n$ such that $\lambda_1+\cdots+\lambda_n \geq 0$
and $S(\lambda)=S$. Here each $\lambda_i=-1$ contributes a
  factor of $x$, each $\lambda_i=1$ contributes a factor of $-1/x$, and
  each $i\in S(\lambda)$ contributes an additional factor of $-1$ to the
  weight of $\lambda$. Just like in \cite[Section
  7.4]{Bayer-Ehrenborg} and in Section~\ref{sec:t-cd} above, we may
associate to each  
$\lambda\in \{-1,1\}^n$ a lattice path starting at $(0,0)$,
such that each $\lambda_i$ is replaced with a step $(1,\lambda_i)$.
We may then use this 
observation to express the polynomials $\st_h(S,x)$ in the basis 
$\{\QQ_k(x)\}_{k\geq 0}$. For this purpose
we need to recall the notion of {\em the unique
sparse interval system ${\mathcal I}[S]$ associated to a set of
positive integers $S$}, see~\cite{Billera-Hetyei-planar}. An antichain of
intervals ${\mathcal I}=\{[i_1,j_1],\ldots,[i_r,j_r]\}$ satisfying
$i_1< \cdots < i_r$ is {\em sparse} if for all $k\in [1,r-1]$ we
have $j_k+1<i_{k+1}$. For every set $S$ of positive integers there is 
a unique sparse interval system ${\mathcal I}$ such that $S$ is the
union of the intervals of ${\mathcal I}$. We denote this family of
intervals by ${\mathcal I}[S]$. We obtain the following short toric
generalization of~\cite[Theorem 3.1]{Bayer-Ehrenborg}.
\begin{proposition}
\label{prop:Fineh-t-exp}
Let $S\subseteq [1,n]$ be a set given by ${\mathcal I}[S]=\{[i_1,j_1],
\cdots,[i_r,j_r]\}$. Then the polynomial $\st_h(S,x)$ appearing in
\eqref{eq:st-h} is given by   
$$
\st_h(S,x)=
\left\{
\begin{array}{ll}
\displaystyle (-1)^{r+\frac{n}{2}} \prod_{k=1}^r \left(C_{\frac{j_k-i_k}{2}}
C_{\frac{i_k-j_{k-1}}{2}}\right) C_{\frac{n-j_r}{2}}
&\mbox{if $j_r<n$;}\\
\displaystyle (-1)^{r+\frac{i_r-1}{2}} \prod_{k=1}^{r-1} C_{\frac{j_k-i_k}{2}}
\prod_{k=1}^r C_{\frac{i_k-j_{k-1}}{2}} \QQ_{n-i_r+1}(x)
&\mbox{if $j_r=n$}.\\
\end{array}
\right.
$$
Here we set $j_0:=0$.
\end{proposition} 
Note that $n$ must be even
if $j_r<n$ since $\lambda_1+\cdots+\lambda_n=0$. Similarly, $i_r$ must
be odd if $j_r=n$.

\end{document}